\documentclass{article}
\usepackage{amsmath,amssymb,amsthm,amscd,latexsym,}
 
\usepackage{graphicx}

\makeatletter
\renewcommand{\mod}[1]{\allowbreak \if@display \mkern 8mu \else
\mkern 5mu\fi {\operator@font mod}\,\,#1}
\makeatother

 \newcommand{\bq}{\mathbb Q}
\newcommand{\br}{\mathbb R}
\newcommand{\bz}{\mathbb Z}
\newcommand{\bff}{\mathbb F}

\newcommand{\bk}{\mathbb K}
\newcommand{\bo}{\mathbb O}

\DeclareMathOperator{\ch}{ch}

\DeclareMathOperator{\discr}{{\rm discr}\,}
\newtheorem{theorem}{Theorem}
\newtheorem{proposition}[theorem]{Proposition}
\newtheorem{definition}[theorem]{Definition}

\newtheorem{lemma}[theorem]{Lemma}

\newcommand\F{\mathcal F}
\newcommand\Hh{\mathcal H}
\newcommand\La{\mathcal L}

\newcommand\M{\mathcal M}

\begin{document}
\title{The transition constant for arithmetic hyperbolic reflection groups}
\date{}
\author{Viacheslav V. Nikulin\footnote{This paper was written with the
financial support of EPSRC, United Kingdom (grant no. EP/D061997/1)}}
\maketitle

\begin{abstract} The transition constant was introduced in our 1981 paper
and denoted as $N(14)$. It is equal to the maximal degree of the ground
fields of V-arithmetic connected edge graphs with 4 vertices and
of the minimality $14$. This constant is fundamental since if the degree
of the ground field of an arithmetic hyperbolic reflection group
is greater than $N(14)$, then the field comes from special plane
reflection groups. In \cite{Nik6}, we claimed its upper bound $56$.
Using similar but more difficult considerations, here we show that
the upper bound is $25$.

As applications, using this result and our methods, we show that the
degree of ground fields of arithmetic hyperbolic reflection groups
in dimensions at least $6$ has the
upper bound  $25$; in dimensions $3,\,4,\,5$ it has the upper bound
$44$. This significantly improves our results in \cite{Nik6, Nik7, Nik8}.
Additionally using recent results by Belolipetsky and Maclachlan,
the last upper bound can be improved to $35$.

In Appendix, we give a review and corrections to Section 1 of our papers
\cite{Nik2} which is important for our methods.
\end{abstract}

\section{Introduction} \label{introduction} The transition constant
was introduced in our paper \cite{Nik2} and was denoted in \cite{Nik2}
as $N(14)$. It is equal to the maximal degree of the ground
fields of V-arithmetic connected edge graphs with 4 vertices and
of the minimality $14$. We showed in \cite{Nik2} that the number
of these fields is finite, and $N(14)$ is a finite effective constant.

It was shown in \cite{Nik2} that the degree of ground fields of
arithmetic hyperbolic reflection groups in dimensions at least 10
has the upper bound $N(14)$. In \cite{Nik6}, we had shown  that
the degree of ground fields of arithmetic  hyperbolic reflection groups
in dimensions at least 6 is bounded by the maximum from $N(14)$ and $11$
(here $11$ is the upper bound for the degree of ground fields
of plain arithmetic hyperbolic reflection groups with quadrangle
fundamental polygon of the minimality $14$).
In \cite{Nik5}, we showed that the degree of ground fields of arithmetic
hyperbolic reflection groups is bounded by the maximum of $N(14)$ and the
degree of ground fields of arithmetic hyperbolic reflection groups in
dimensions two and three.
In general, this constant is fundamental since if the degree
of the ground field of an arithmetic hyperbolic reflection group
is greater than $N(14)$, then the field comes from special plane
reflection groups. See \cite{Nik2}, \cite{Nik5} --- \cite{Nik8}.

In \cite{Nik6}, we claimed that $N(14)\le 56$. Here we improve
this bound. Using similar but more difficult considerations,
we show that $N(14)\le 25$.

As applications (see Section \ref{Applications}), using our methods,
we show that the degree of ground fields of arithmetic hyperbolic
reflection groups in dimensions $n\ge 6$ has the upper bound $25$;
in dimensions $n=3,\,4,\,5$, it has the upper bound $44$.
Remark that this also gives another proof of finiteness in dimension
$n=3$ which is different from the first proof by I. Agol \cite{Agol}.
These significantly improve results of our recent papers
\cite{Nik6}---\cite{Nik8}.

In \cite{B}, Belolipetsky obtained the upper bound $35$ for $n=3$.
In \cite{M}, Maclachlan obtained the upper bound $11$ for $n=2$.
Using these results and our results from \cite{Nik5}, we can improve
the upper bound $44$ in dimensions $n=4,\,5$. We show that the
upper bound is $35$ in these dimensions.

We hope that our results will be important for further classification.
See an example in Theorem \ref{class1}.

In Appendix, we review and correct results of Section 1 of
our paper \cite{Nik2} which are very important for our methods.

I am grateful to Professor C. Maclachlan for sending me his recent
preprint \cite{M}.

\section{Reminding of some basic facts about hyperbolic fundamental
polyhedra}\label{sec:V-arpol}

Here we remind some basic definitions and results about fundamental
chambers (always for discrete reflection groups) in hyperbolic spaces and
their Gram matrices.
See \cite{Vin1}, \cite{Vin5} and \cite{Nik1}, \cite{Nik2}.

We work with Klein model of a hyperbolic space
$\La$ associated to a hyperbolic form $\Phi$ over the field of
real numbers $\br$ with signature $(1,n)$, where $n=\dim \La$.
Let $V=\{x\in \Phi|x^2>0\}$ be the cone determined by $\Phi$,
and let $V^+$ be one of the two halves of this cone. Then
$\La=\La(\Phi)=V^+/\br^+$ is the set of rays in $V^+$; we let
$[x]$ denote the element of $\La$ determined by the ray $\br^+x$
where $x\in V^+$ and $\br^+$ is the set of all positive real numbers.
The hyperbolic distance is given by the formula
$$
\rho([x],[y])=(x\cdot y)/\sqrt{x^2y^2},\ [x], [y]\in \La,
$$
then the curvature of $\La$ is equal to $-1$.

Every half-space $\Hh^+$ in $\La$
determines and is determined by the orthogonal element $e\in \Phi$ with
square $e^2=-2$:
$$
\Hh^+=\Hh_e^+=\{[x]\in \La|x\cdot e\ge 0\}.
$$
It is bounded by the hyperplane
$$
\Hh^+=\Hh_e^+=\{[x]\in \La|x\cdot e = 0\}
$$
orthogonal to $e$.  If two half-spaces
$\Hh_{e_1}^+$, $\Hh_{e_2}^+$ where $e_1^2=e_2^2=-2$
have a common non-empty open subset in $\La$, then
$\Hh_{e_1}\cap \Hh_{e_2}$ is an angle of the value $\phi$ where
$2\cos{\phi}=e_1\cdot e_2$ if $-2<e_1\cdot e_2\le 2$,
and the distance between hyperplanes $\Hh_{e_1}$ and
$\Hh_{e_2}$ is equal to $\rho$ where $2\ch{\rho}=e_1\cdot e_2$
if $e_1\cdot e_2>2$.

A convex polyhedron $\M$ in $\La$ is intersection of a finite
number of half-spaces $\Hh^+_e$, $e\in P(\M)$, where $P(\M)$ are all the
vectors with square $-2$ which are orthogonal to the faces
(of the codimension one) of $\M$ and are directed outward. The matrix
\begin{equation}
A=(a_{ij})=(e_i\cdot e_j),\ e_i, e_j\in P(\M),
\label{Grammatr}
\end{equation}
is the Gram matrix $\Gamma (\M)=\Gamma (P(\M))$ of $\M$.
It determines $\M$ uniquely up to motions
of $\La$. If $\M$ is sufficiently general, then $P(\M)$ generates $\Phi$,
and the form $\Phi$ is
\begin{equation}
\Phi=\sum_{e_i,e_j\in P(\M)}{a_{ij}X_iY_j}\mod {\rm Kernel},
\label{formPhi}
\end{equation}
and $P(\M)$ naturally identifies with a subset of $\Phi$ and defines $\M$.

The polyhedron $\M$ is a fundamental chamber of a discrete reflection group
$W$ in $\La$ if and only if $a_{ij}\ge 0$ and
$a_{ij}=2\cos {\frac{\pi}{m_{ij}}}$ where $m_{ij}\ge 2$ is an integer if
$a_{ij}<2$ for all $i\not=j$.
Symmetric real matrices $A$ satisfying these conditions and having all their
diagonal elements equal to $-2$ are called {\it fundamental} (then the set
$P(\M)$ formally corresponds to indices of the matrix $A$).
As usual, further we identify  fundamental matrices with fundamental
graphs $\Gamma$. Their vertices correspond to
$P(\M)$. Two different vertices $e_i\not=e_j\in P(\M)$ are connected by
the thin edge of the integer weight $m_{ij}\ge 3$ if
$0<a_{ij}=2\cos {\frac{\pi}{m_{ij}}}<2$, by the thick edge if $a_{ij}=2$,
and by the broken edge of the weight $a_{ij}$ if $a_{ij}>2$.
In particular, the vertices $e_i$ and $e_j$ are disjoint if and only if
$e_i\cdot e_j=a_{ij}=2\cos{\frac{\pi}{2}}=0$. Equivalently, $e_i$ and
$e_j$ are perpendicular (or orthogonal). See some examples of such graphs
in Figures \ref{figLanner} --- \ref{graphg45} below.

For a real $t>0$, we say that a fundamental matrix $A=(a_{ij})$ (and the
corresponding fundamental chamber $\M$) {\it has minimality $t$}
if $a_{ij}<t$ for all $a_{ij}$. Here we follow \cite{Nik1},
\cite{Nik2}. Further, the minimality $t=14$ will be especially important.

It is known that fundamental domains of arithmetic hyperbolic
groups must have finite volume. Let us assume that it is valid for a
fundamental chamber $\M$ of a hyperbolic discrete reflection group.
As Vinberg had shown \cite{Vin1}, in order for $\M$ to be a fundamental
chamber of an arithmetic reflection group $W$ in $\La$, it is necessary and
sufficient that
all of the cyclic products
\begin{equation}
b_{i_1\dots i_m}=a_{i_1i_2}\cdot a_{i_2i_3}\cdots a_{i_{m-1}i_m}\cdot
a_{i_mi_1}
\label{cycprod}
\end{equation}
be algebraic integers, that the field
$\widetilde{\bk}=\bq(\{a_{ij}\})$ be
totally real, and that, for any embedding $\widetilde{\bk}\to \br$ not
the identity over the {\it ground field} $\bk=\bq(\{b_{i_1\dots i_m}\})$
generated by all of the cyclic products \eqref{cycprod}, the form
\eqref{formPhi}
be negative definite.

Fundamental real matrices $A=(a_{ij})$, $a_{ij}=e_i\cdot e_j$,
$e_i,\,e_j\in P(\M)$ (or the corresponding graphs),
with a hyperbolic form $\Phi$ in \eqref{formPhi} and satisfying
these Vinberg's conditions will be further called
{\it V-arithmetic} (here we don't require that the corresponding
hyperbolic polyhedron $\M$ has finite volume). It is well-known
(and easy to see; see arguments in Sect. \ref{subsec:gam(4)1}) that
a subset $P\subset P(\M)$ also defines a V-arithmetic
matrix $(e_i\cdot e_j)$, $e_i, e_j\in P$, with the same ground
field $\bk$ if the subset $P$ is hyperbolic, i. e. the corresponding to
$P$ form \eqref{formPhi} is hyperbolic.

\section{V-arithmetic edge polyhedra}{\label{V-arithedgepol}

A fundamental chamber $\M$ (and the corresponding Gram matrix $A$ or
a graph) is called {\it edge chamber (matrix, graph)} if all hyperplanes
$\Hh_e$, $e\in P(\M)$, contain one of two distinct vertices $v_1$ and
$v_2$ of the 1-dimensional edge $v_1v_2$ of $\M$. Assume that both vertices
$v_1$ and $v_2$ are finite (further we always consider this case).
Further we call this edge chambers {\it finite}. Assume that
$\dim \La=n$. Then $P(\M)$ consists of
$n+1$ elements: $e_1$, $e_2$ and $n-1$ elements $P(\M)-\{e_1,e_2\}$.
Here $P(\M)-\{e_1,e_2\}$ corresponds to hyperplanes which contain the edge
$v_1v_2$ of $\M$. The $e_1$ corresponds to the hyperplane which contains
$v_1$ and does not contain $v_2$. The $e_2$ corresponds to the hyperplane
which contains $v_2$ and does not contain $v_1$. Then the set $P(\M)$
is hyperbolic (it has hyperbolic Gram matrix),
but its subsets $P(\M)-\{e_1\}$ and $P(\M)-\{e_2\}$ are negative definite
(they have negative definite Gram matrix) and define Coxeter graphs.
Only the element $u=e_1\cdot e_2$ of the Gram matrix of $\M$ can be greater
than $2$. Thus, $\M$ will have the minimality $t$ if and only if
$u=e_1\cdot e_2<t$.

From considerations above, the Gram graph $\Gamma (P(\M))$ of an edge chamber
has only one hyperbolic connected component $P(\M)^{hyp}$ (containing
$e_1$ and $e_2$) and several negative definite connected components.
Gram matrix $\Gamma (P(\M)^{hyp})$ evidently also corresponds to an
edge chamber of the dimension $\#P(\M)^{hyp}-1$. If  $\M$ is
V-arithmetic, the ground field $\bk$ of $\M$ is the same as for
the hyperbolic connected component $\Gamma (P(\M)^{hyp})$

The following result had been proved in \cite{Nik2}.

\begin{theorem} (\cite[Theorem 2.3.1]{Nik2}) Given $t>0$, there exists an
effective constant $N(t)$ such that every
V-arithmetic edge chamber of the minimality $t$ with
ground field $\bk$ of degree greater that $N(t)$ over $\bq$ has the hyperbolic
connected component of its Gram graph
which has less than $4$ elements.
\label{tholdVedge}
\end{theorem}

Considerations in \cite{Nik2} (and also \cite{Nik1})
also show that the set of possible
ground fields $\bk$ of hyperbolic
connected components with at least $4$ vertices
of V-arithmetic edge chambers of minimality $t$ is finite.
Even the set of Gram graphs
$\Gamma(P(\M)^{hyp})$ of minimality $t$
with fixed $\ge 4$ number of vertices is finite.
Taking this under consideration, here we want to formulate and prove
more efficient variant of this theorem. We restrict by the minimality
$t=14$ to get an exact estimate for the constant $N(14)$, but the
same finiteness results are valid for any $t>0$.

Following \cite{Nik6}, below we formulate a more efficient definition of the
constants $N(t)$, $t>0$, and $N(14)$.

\begin{figure}
\begin{center}
\includegraphics[width=6cm]{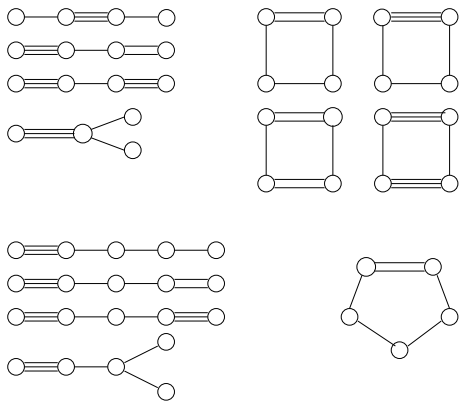}
\end{center}
\caption{All arithmetic Lann\'er graphs with at least $4$ vertices}
\label{figLanner}
\end{figure}

\subsection{Arithmetic Lann\'er graphs with $\ge 4$ elements}
\label{subsecLanner}

We remind that Lann\'er graphs are Gram graphs of bounded
fundamental hyperbolic simplexes. They are characterized as
hyperbolic fundamental graphs such that any their proper subgraph
is a Coxeter graph. They were classified by Lann\'er \cite{Lan}.
In Figure \ref{figLanner} we give all arithmetic Lann\'er graphs
with at least $4$ vertices (only one Lann\'er graph with $\ge 4$
vertices is not arithmetic). As usual, we replace a thin edge of
the weight $k$ by $k-2$-edges for a small $k$. Ground fields of
Lann\'er graphs with $\ge 4$ vertices give three fields:
\begin{equation}
\F L^4=\{\bq,\,\bq(\sqrt{2}),\ \bq(\sqrt{5})\}.
\label{fLanner4}
\end{equation}
See \cite{Vin3} for details.

\subsection{Arithmetic triangle graphs}
\label{subsectriangle}

\begin{figure}
\begin{center}
\includegraphics[width=6cm]{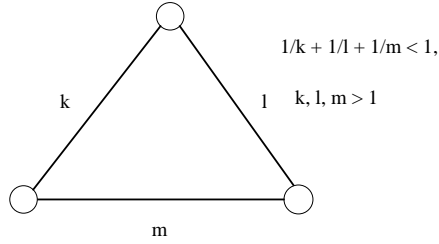}
\end{center}
\caption{Triangle graphs}
\label{figtriangles}
\end{figure}

Triangle graphs are Gram graphs of bounded fundamental triangles on
hyperbolic plane (we don't consider non-bounded triangles).
Equivalently, they are Lann\'er graphs with $3$
vertices. They are given in Figure \ref{figtriangles} where
$2\le k,l,m$ and
$$
\frac{1}{k}+\frac{1}{l}+\frac{1}{m}<1.
$$
Arithmetic triangles were enumerated
by Takeuchi \cite{Tak2}. All bounded arithmetic triangles
are given by the following triplets $(k,l,m)$:
$$
(2,3,7\ -\ 12),\ (2,3,14),\ (2,3,16),\ (2,3,18),\ (2,3,24),\ (2,3,30),\
(2,4,5\ -\ 8),
$$
$$
(2,4,10),\ (2,4,12),\ (2,4,18),\ (2,5,5),\ (2,5,6),\ (2,5,8),\
(2,5,10),\ (2,5,20),
$$
$$
(2,5,30),\ (2,6,6),\ (2,6,8),\ (2,6,12),\ (2,7,7),\ (2,7,14),\
(2,8,8),\
(2,8,16),
$$
$$
(2,9,18),\ (2,10,10),\ (2,12,12),\ (2,12,24),\ (2,15,30),\ (2,18,18),\
(3,3,4\ -\ 9),\
$$
$$
(3,3,12),\ (3,3,15),\ (3,4,4),\ (3,4,6),\ (3,4,12),\ (3,5,5),\ (3,6,6),\
(3,6,18),
$$
$$
(3,8,8),\ (3,8,24),\ (3,10,30),\ (3,12,12),\ (4,4,4\ - 6),\ (4,4,9),\
(4,5,5),\ (4,6,6),
$$
$$
(4,8,8),\ (4,16,16),\ (5,5,5),\ (5,5,10),\ (5,5,15),\ (5,10,10),\
(6,6,6),\ (6,12,12),
$$
$$
(6,24,24),\ (7,7,7),\ (8,8,8),\ (9,9,9),\ (9,18,18),\ (12,12,12),\ (15,15,15).
$$
Their ground fields were found by Takeuchi \cite{Tak3}. They give the set
of fields
\begin{equation}
\begin{array}{l}
\F T=
\{\bq\}\cup \{\bq(\sqrt{a})\ |\ a=2,\,3,\,5,\,6\}\cup
\{\bq(\sqrt{2},\sqrt{3}),\ \bq(\sqrt{2},\sqrt{5})\}\cup\\
\cup \{\bq(\cos{\frac{2\pi}{b}})\ |\ b=7,\,9,\,11,\,15,\,16,\,20\}.
\end{array}
\label{ftriangles}
\end{equation}

\subsection{V-arithmetic connected finite edge graphs with $4$ vertices
for $2<u<14$}\label{subsecedggraphs4}

Using classification of Coxeter graphs, it is easy to draw
all possible pictures of connected finite edge graphs $\Gamma^{(4)}$
with $4$ vertices and $u=e_1\cdot e_2>2$. They correspond to all
3-dimensional finite fundamental edge polyhedra with connected Gram
graph and $u>2$. They are given in Figure \ref{figfivegraphs} and give
five types of graphs $\Gamma=\Gamma^{(4)}_i$, $i=1,\,2,\,3,\,4,\,5$.
All possible natural parameters $s,\,k,\,r,\,p\ge 2$ for these graphs can
be easily enumerated by the condition that
$\Gamma-\{ e_1 \}$, $\Gamma-\{ e_2\}$ are Coxeter graphs.
They will be given in Sec. \ref{secfieldedggraphs4}  below.

\begin{definition}
For $i=1,\,2,\,3,\,4,\,5$ and $t>0$ we denote by $\Gamma^{(4)}_i(t)$ the set of
all V-arithmetic connected finite
edge graphs with 4
vertices  $\Gamma^{(4)}_i$ of the minimality $t$, i. e. for $2<u<t$, and by
$$
\F \Gamma^{(4)}_i(t)
$$
the set of all their ground fields.
\label{fieldsedge4}
\end{definition}

All V-arithmetic graphs $\Gamma^{(4)}_i$ for $2<u<t$  give particular
cases of graphs of V-arithmetic edge polyhedra
with hyperbolic connected component having $4$ vertices and minimality $t$.
Thus, by Theorem \ref{tholdVedge}, degree (over $\bq$) of fields from
$\F\Gamma^{(4)}_i(t)$ is bounded by the effective
constant $N(t)$. It follows that the sets of V-arithmetic graphs
$\Gamma^{(4)}_i(t)$ and fields $\F \Gamma^{(4)}_i(t)$ are
also finite.

Vice verse, Theorem \ref{fieldsedge4} can be deduced from finiteness
of the sets of fields above because of the
following easy statement.

\begin{figure}
\begin{center}
\includegraphics[width=10cm]{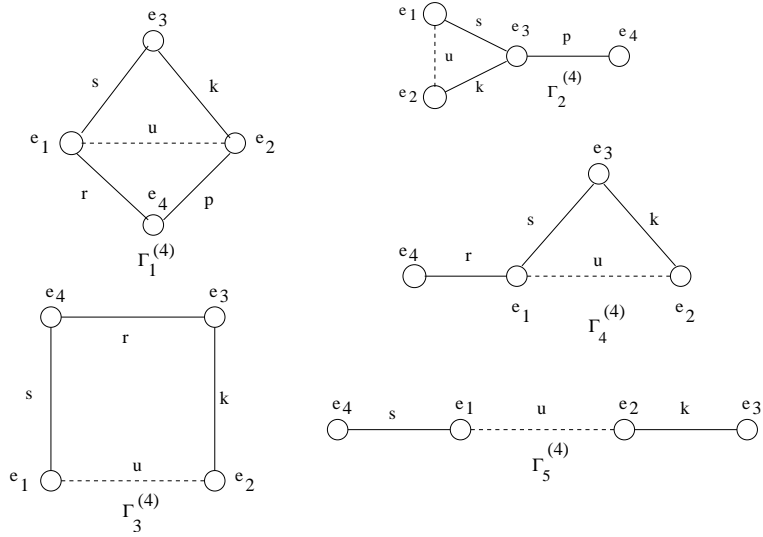}
\end{center}
\caption{Five graphs $\Gamma^{(4)}_i$, $i=1,\,2,\,3,\,4,\,5$.}
\label{figfivegraphs}
\end{figure}

\begin{proposition} The ground field of any
V-arithmetic edge chamber of the minimality $t>0$ with the hyperbolic
connected component of its Gram graph  having at least $4$ vertices
belongs to one of the finite sets of
fields $\F L^4$, $\F T$ and $\F\Gamma^{(4)}_i(t)$, $1\le i\le 5$,
introduced above.

In particular, Theorem \ref{tholdVedge} is equivalent to finiteness of
the sets of fields
$\F\Gamma^{(4)}_i(t)$, $ i=1,\,2,\,3,\,4,\,5$.
\label{prnewVedge}
\end{proposition}

\begin{proof} See \cite{Nik6}.
\end{proof}

Degree of fields from $\F L^4$ is  bounded by $2$, and degree of fields from
$\F T$ is bounded by $5$.

For arithmetic hyperbolic reflection groups, the minimality $t=14$ is
especially important.
Using the same methods as for the proof of Theorem \ref{tholdVedge}
in \cite{Nik2},  we can prove
the following effective upper bounds which improve the upper bounds
which we claimed in \cite{Nik6} (e.g., in \cite{Nik6} we claimed
that $N(14)\le 56$).

\begin{theorem} The degree of fields from $\F \Gamma^{(4)}_1(14)$ is
bounded ($\le$) by $23$.

The degree of fields from $\F \Gamma^{(4)}_2(14)$ is bounded by $23$.

The degree of fields from $\F \Gamma^{(4)}_3(14)$ is bounded by $25$.

The degree of fields from $\F \Gamma^{(4)}_4(14)$ is bounded by $25$.

The degree of fields from $\F \Gamma^{(4)}_5(14)$ is bounded by $25$.

Thus, the constant $N(14)$ of Theorem \ref{tholdVedge} can be taken to
be $N(14)=25$.
\label{thdegrees}
\end{theorem}

\begin{proof} The proof will be given in the next Section
\ref{secfieldedggraphs4}.
\end{proof}

\section{Ground fields of V-arithmetic connected finite edge graphs with
four vertices  of the minimality $14$.}
\label{secfieldedggraphs4}

Here we shall obtain explicit upper bounds of  degrees of fields from
the finite sets $\F \Gamma^{(4)}_i(14)$, $1\le i\le 5$ (see Definition
\ref{fieldsedge4}), and prove Theorem \ref{thdegrees}.  Moreover, our
considerations will deliver important information about these
sets of fields.

\subsection{Some results on hyperbolic numbers}
\label{subsec:hypnumbers}

Like for the proof of Theorem \ref{tholdVedge} from \cite{Nik2}, we use
the following general results from \cite[Section 1]{Nik2}
(see Section \ref{sec:Appendix} which also contains some corrections).

\begin{theorem} (\cite[Theorem 1.2.1]{Nik2}) Let $\bff$ be a totally
real algebraic number field, and let each
embedding $\sigma:\bff\to \br$ corresponds to an interval
$[a_\sigma,b_\sigma]$ in $\br$ where
$$
\prod_{\sigma }{\frac{b_\sigma-a_\sigma}{4}}<1.
$$
In addition, let the natural number $m$ and the intervals
$[s_1,t_1],\dots, [s_m,t_m]$ in $\br$ be fixed. Then there exists
a constant $N(s_i,t_i)$ such that, if $\alpha$ is a totally real
algebraic integer and if the following inequalities hold for the
embeddings $\tau:\bff(\alpha) \to \br$:
$$
s_i\le \tau(\alpha)\le t_i\ \ for\ \ \tau=\tau_1,\dots ,\tau_m,
$$
$$
a_{\tau | \bff}\le \tau(\alpha)\le b_{\tau | \bff}\ \ for\ \
\tau\not=\tau_1,\dots,\tau_m,
$$
then
$$
[\bff(\alpha):\bff]\le N(s_i,t_i).
$$
\label{th121}
\end{theorem}

\begin{theorem} (\cite[Theorem 1.2.2]{Nik2})
Under the conditions of Theorem \ref{th121}, \newline $N(s_i,t_i)$ can be
taken to be $N(s_i,t_i)=N$,
 where $N$ is the least natural number
solution of the inequality
\begin{equation}
N\ln{(1/R)} - M\ln{(2N+2)}-\ln{B}\ge \ln{S}.
\label{cond for n}
\end{equation}

Here
\begin{equation}
M=[\bff : \bq],\ \ \ B=\sqrt{|{\rm discr\ } \bff|};
\label{MB}
\end{equation}
\begin{equation}
R=\sqrt{\prod_\sigma {\frac{b_\sigma-a_\sigma}{4}}},\ \ \
S=\prod_{i=1}^{m}{\frac{2er_i}{b_{\sigma_i}-a_{\sigma_i}}}
\label{RS}
\end{equation}
where
\begin{equation}
\sigma_i=\tau_i|\bff,\ \ \  r_i=\max\{{|t_i-a_{\sigma_i}|,
|b_{\sigma_i}-s_i|}\}.
\label{ri}
\end{equation}
\label{th122}
\end{theorem}

We note that the proof of Theorems \ref{th121} and \ref{th122}
uses a variant of Fekete's Theorem (1923) about existence of
non-zero integer polynomials of bounded degree which differ only
slightly from zero on appropriate intervals. See \cite[Theorem
1.1.1]{Nik2} (see its corrections in Section \ref{sec:Appendix}, Theorems
\ref{thFekete1}, \ref{thFekete2}).

Unfortunately, Theorems \ref{th121} and \ref{th122} usually
give a poor upper bound for the degree.
They should mainly be considered as existence theorems.
Usually we shall use explicit polynomials to bound
the degree in similar cases.

We use the following statement.
Similar arguments one should use to prove Theorems
\ref{th121} and \ref{th122}, see Section \ref{sec:Appendix}.

\begin{lemma}
Let $\bff$ be a totally
real algebraic number field, and let each
embedding $\sigma:\bff\to \br$ corresponds to an interval
$[a_\sigma,b_\sigma]$ in $\br$.
In addition, let the natural number $m$ and the intervals
$[s_1,t_1],\dots, [s_m,t_m]$ in $\br$ be fixed.

Let $P(x)$ be a non-zero polynomial over the ring of
integers of $\bff$, and
$$
\delta(\sigma)=\max_{x\in [a_\sigma,b_\sigma]}{|P^\sigma(x)|}\ \ \text{for}\ \
\sigma:\bff\to \br,
$$
and
$$
a_i=\max_{x\in [s_i,t_i]}{|P^{(\sigma_i|\bff)}(x)|}.
$$
Suppose that $\prod_\sigma {\delta(\sigma)}<1$.

Let $\alpha$ be a totally real algebraic integer such that the
following inequalities hold for all
embeddings $\tau:\bff(\alpha) \to \br$:
$$
s_i\le \tau(\alpha)\le t_i\ \ for\ \ \tau=\tau_1,\dots ,\tau_m,
$$
$$
a_{\tau | \bff}\le \tau(\alpha)\le b_{\tau | \bff}\ \ for\ \
\tau\not=\tau_1,\dots,\tau_m,
$$
Then we have the following bound for $[\bff(\alpha):\bff]$:
$$
1\le [\bff(\alpha):\bff]\le
\frac{\ln{\prod_i{a_i}}-\ln{\prod_i{\delta(\tau_i|\bff)}}}
{-\ln{\prod_\sigma {\delta(\sigma)}}}\
$$
if $P(\alpha)\not=0$. Note that $\alpha$ does not exist if the
right hand side is $<1$.
\label{lempolmain}
\end{lemma}

\begin{proof} Since $P(\alpha)\not=0$, we have the inequalities for
the norm
$$
1\le |N_{\bff(\alpha)/\bq}(P(\alpha))|=
$$
$$
|\prod_\tau{\tau(P(\alpha))|}=\prod_{\tau}{|P^\tau(\tau(\alpha))|}=
\prod_{\tau\not=\tau_i}{P^\tau(\tau(\alpha))}
\prod_{i=1}^{m}{|P^{\tau_i}(\tau_i(\alpha))|}
$$
$$
\le \prod_{\tau\not=\tau_i}{\max_{[a_{\tau|\bff},b_{\tau|\bff}]}{|P^\tau(x)|}}
\prod_{i=1}^m{\max_{[s_i,t_i]}{|P^{\sigma_i}(x)|}}
$$
$$
\le \prod_{\tau}{\max_{[a_{\tau|\bff},b_{\tau|\bff}]}{|P^\tau(x)|}}
\prod_{i=1}^m{\max_{[s_i,t_i]}{|P^{\sigma_i}(x)|}}/
\prod_{i=1}^{i=m}{\max_{[a_{\tau_i|\bff},b_{\tau_i|\bff}]}{|P^{\tau_i}(x)|}}=
$$
$$
\le \left(\prod_{\sigma}{\delta(\sigma)}\right)^{[\bff(\alpha):\bff]}
\prod_{i=1}^m{a_i}/
\prod_{i=1}^{i=m}{\delta(\tau_i)}.
$$
Taking logarithm from both sides and using that
$\prod_\sigma {\delta(\sigma)}<1$, we obtain the statement.

This is similar to the proof of Theorems
\ref{th121} and \ref{th122} in Section \ref{sec:Appendix}.
\end{proof}

\subsection{Ground fields of some V-arithmetic connected edge graphs
with 3 vertices with the given minimality.}
\label{subsec:gam(3)}

In spite of we are mainly interested in connected edge
graphs with four vertices, some connected edge graphs with three
vertices will be very important.

All connected edge graphs with three vertices are given in Figure
\ref{graphg3}. They are $\Gamma_1^{(3)}$ where $s, k\ge 3$ and
$\Gamma_2^{(3)}$ where $d\ge 3$. We denote by $\F \Gamma_1^{(3)}(t)$
and $\F \Gamma_2^{(3)}(t)$ their ground fields for the minimality
$t>2$. It means that $2<\sigma^{(+)}(u)<t$.

\subsubsection{Ground fields of $\Gamma_2^{(3)}(t)$}
\label{subsubsec:Gamma2(3)}
First let us consider the graphs $\Gamma_2^{(3)}(t)$.
The corresponding Gram matrix is
\begin{equation}
\left(\begin{array}{ccc}
-2     & u                    &   0                   \\
 u     & -2                   & 2\cos{\frac{\pi}{d}}  \\
 0     &  2\cos{\frac{\pi}{d}}  & -2
\end{array}\right)
\label{gramm2(3)}
\end{equation}
where $d\ge 3$ is an integer, $u$ is a totally real algebraic integer.
The ground field $\bk=\bq(u^2,\cos^2{(\pi/d)})$.
The determinant $d(u)$ of the Gram matrix is given by the equality
$$
\frac{d(u)}{2}=u^2-4\sin^2{(\pi/d)}.
$$
It follows that $\Gamma_2^{(3)}$ is V-arithmetic of the minimality $t>2$,
if and only if for $\alpha=u^2$ we have
\begin{equation}
0<\sigma(\alpha)<\sigma(4\sin^2{\frac{\pi}{d}})
\label{F32cond1}
\end{equation}
for all $\sigma:\bk\to \br$ which are different from $\sigma^{(+)}$, and
\begin{equation}
4<\sigma^{(+)}(\alpha)<t^2.
\label{F32cond2}
\end{equation}
It follows that $\bff_d=\bq(\cos^2{(\pi/d)})\subset \bk=\bq(\alpha)$.
We will be especially interested in $t=14$ and $t=16$.
Let us estimate the degree $n=[\bk:\bq]$.


\begin{figure}
\begin{center}
\includegraphics[width=6cm]{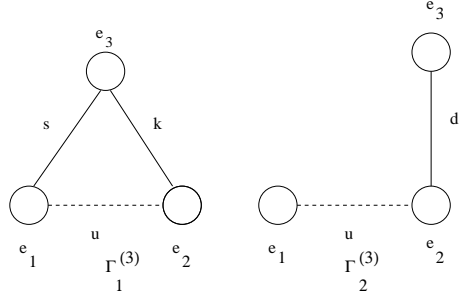}
\end{center}
\caption{Connected edge graphs with 3 vertices}
\label{graphg3}
\end{figure}

Assume that $d=3$. Then $\alpha$ satisfies
$$
0<\sigma(\alpha)<3,\ \ \ \ \ 4<\sigma^{(+)}(\alpha)<t^2.
$$
We take the polynomial
\begin{equation}
P(x)=x^3(x-1)^4(x-2)^4(x-3)^3(x^2-3x+1)^3.
\label{pol3}
\end{equation}
of the degree $20$. The maximum of $|P(x)|$ on the interval $[0,3]$
is equal to $\delta=884736/9765625=0.0905969664$. Using $P(x)$,
by Lemma \ref{lempolmain}, we get the upper bound
$$
[\bk:\bq]\le 1+\frac{\ln(P(t^2))}{(-\ln{(\delta)})}.
$$
For $t=14$, we get $[\bk:\bq]\le 44$, and for $t=16$,
we get $[\bk:\bq]\le 47$. (Theorems \ref{th121} and \ref{cond for n}
give more poor results: $76$ for $t=14$ and $78$ for $t=16$.)

Assume that $d=4$. Then $\alpha$ satisfies the inequalities
$$
0<\sigma(\alpha)<2,\ \ \ \ \ 4<\sigma^{(+)}(\alpha)<t^2.
$$
We take the polynomial
\begin{equation}
P(x)=x(x-1)^2(x-2).
\label{pol4}
\end{equation}
of the degree $4$. The maximum of $|P(x)|$ on the interval $[0,2]$
is equal to $\delta=1/4$. Using $P(x)$, by Lemma \ref{lempolmain},
we obtain the bound
$$
[\bk:\bq]\le 1+\frac{\ln(P(t^2))}{(-\ln(\delta))}.
$$
For $t=14$ and $t=16$, we get $[\bk:\bq]\le 16$.
(Theorems \ref{th121} and \ref{cond for n} give more poor bound:
$31$.)

Assume that $d=5$.
We take the polynomial
\begin{equation}
P(x)=x(x-1)\left(x+1-4\sin^2{\frac{\pi}{5}}\right)
\left(x-4\sin^2{\frac{\pi}{5}}\right).
\label{pol5}
\end{equation}
of the degree $4$ over the ring of integers of $\bff_5$.
The maximum of $|P(x)|$ on the
interval $[0,4\sin^2(\pi/5)]$ is equal to $\delta_1=0.04559...$. For its
conjugate polynomial
$$
P^\sigma(x)=x(x-1)\left(x+1-4\sin^2{\frac{2\pi}{5}}\right)
\left(x-4\sin^2{\frac{2\pi}{5}}\right),
$$
the maximum of $|P^\sigma(x)|$
on the conjugate interval $[0,4\sin^2{(2\pi/5)}]$ is equal to
$\delta_2=2.1419...$. Evidently, $P(\alpha)$ is not zero. By
Lemma \ref{lempolmain}, we get the bound
$$
[\bk:\bff_5]\le \frac{\ln{P(t^2)}-\ln{\delta_1}}{-\ln{(\delta_1\cdot \delta_2)}}\ .
$$
For $t=14$ and $t=16$, we get $[\bk:\bff_5]\le 10$, and $[\bk:\bq]\le 20$.
(Theorems \ref{th121} and \ref{cond for n} give the bounds
$[\bk:\bff_5]\le 27$ for $t=14$, and $[\bk:\bff_5]\le 28$ for $t=16$ only.)

Assume that $d\ge 6$. We take the polynomial
\begin{equation}
P(x)=x\left(x-4\sin^2{\frac{\pi}{d}}\right).
\label{polyngen}
\end{equation}
over the ring of integers of  $\bff_d$.
For $\sigma:\bff_d\to \br$, the maximum of $|P^\sigma(x))|$ on the interval
$[0,\sigma(4\sin^2{(\pi/d)})]$ is equal to $(4\sigma(\sin^4{(\pi/d)})$.
Obviously, $P(\alpha)\not=0$.
Then we have (in fact, we repeat the proof of Lemma \ref{lempolmain}
for this case):
$$
1\le |N_{\bk/\bq}(P(\alpha))|<N_{\bff_d/\bq}(4\sin^4{(\pi/d)})^{[\bk:\bff_d]}
\frac{|P(t^2)|}{4\sin^4{(\pi/d)}}=
$$
$$
=\left(\frac{\gamma(d)^2}{4^{[\bff_d:\bq]}}\right)^{[\bk:\bff_d]}
\frac{|P(t^2)|}{4\sin^4{(\pi/d)}}.
$$
where for $d\ge 3$ we have
\begin{equation}
N_{\bff_l/\bq }(4\sin^2{(\pi/d)})=\gamma(d)=\left\{
\begin{array}{cl}
p &\ {\rm if}\  d=p^t>2\   {\rm where\ } p\  {\rm is\  prime,} \\
1 &\ {\rm otherwise,}
\end{array}\
\right .
\label{defgam}
\end{equation}
and $[\bff_d:\bq]=\varphi(d)/2$ where $\varphi(d)$ is the Euler function.

Denoting $m=[\bk:\bff_d]$, we obtain the inequality
$$
m\left(\varphi(d)\ln{2}-2\ln{\gamma(d)}\right)<
\ln{t^2}+\ln{(t^2-4\sin^2{(\pi/d)})}-
\ln({4\sin^4{(\pi/d)})}\ .
$$
It is easy to see that $\varphi(d)\ln{2}-2\ln{\gamma(d)}>0$ for $d\ge 6$.
Thus, we obtain
\begin{equation}
1\le m<\frac{\ln{t^2}+\ln{(t^2-4\sin^2{(\pi/d)})}-\ln({4\sin^4{(\pi/d)})}}
{\varphi(d)\ln{2}-2\ln{\gamma(d)}}
\label{estfor3}
\end{equation}
Respectively $[\bk:\bq]\le m\varphi(d)/2$ where $m$ satisfies
 \eqref{estfor3}.

We use the trivial bound: $\varphi(d)\ge \sqrt{d-2}$. Then
for $t=14$, we obtain that only the following $d\ge 6$ are possible,
and we get the corresponding estimates for $[\bk:\bq]$:

\medskip

\noindent
$d=6$, then $[\bk:\bq]\le 8$;
$d=7$, then $[\bk:\bq]\le 138$;
$d=8$, then $[\bk:\bq]\le 18$;
$d=9$, then $[\bk:\bq]\le 18$;
$d=10$, then $[\bk:\bq]\le 10$;
$d=11$, then $[\bk:\bq]\le 30$;
$d=12$, then $[\bk:\bq]\le 10$;
$d=13$, then $[\bk:\bq]\le 24$;
$d=14$, then $[\bk:\bq]\le 9$;
$d=15$, then $[\bk:\bq]\le 8$;
$d=16$, then $[\bk:\bq]\le 8$;
$d=17$, then $[\bk:\bq]\le 16$;
$d=18$, then $[\bk:\bq]\le 9$;
$d=19$, then $[\bk:\bq]\le 18$;
$d=20$, then $[\bk:\bq]\le 8$;
$d=21$, then $[\bk:\bq]\le 12$;
$d=22$, then $[\bk:\bq]\le 10$;
$d=23$, then $\bk=\bff_{23}$;
$d=24$, then $[\bk:\bq]\le 12$;
$d=25$, then $\bk=\bff_{25}$;
$d=26$, then $[\bk:\bq]\le 12$;
$d=27$, then $\bk=\bff_{27}$;
$d=28$, then $[\bk:\bq]\le 12$;
$d=29$, then $\bk=\bff_{29}$;
$d=30$, then $[\bk:\bq]\le 12$;
$d=31$, then $\bk=\bff_{31}$;
$d=32$, then $\bk=\bff_{32}$;
$d=33$, then $\bk=\bff_{33}$;
$d=34$, then $\bk=\bff_{34}$;
$d=35$, then $\bk=\bff_{35}$;
$d=36$, then $[\bk:\bq]\le 12$;
$d=37$, then $\bk=\bff_{37}$;
$d=38$, then $\bk=\bff_{38}$;
$d=39$, then $\bk=\bff_{39}$;
$d=40$, then $\bk=\bff_{40}$;
$d=42$, then $[\bk:\bq]\le 12$;
$d=44$, then $\bk=\bff_{44}$;
$d=45$, then $\bk=\bff_{45}$;
$d=46$, then $\bk=\bff_{46}$;
$d=48$, then $\bk=\bff_{48}$;
$d=50$, then $\bk=\bff_{50}$;
$d=52$, then $\bk=\bff_{52}$;
$d=54$, then $\bk=\bff_{54}$;
$d=56$, then $\bk=\bff_{56}$;
$d=58$, then $\bk=\bff_{58}$;
$d=60$, then $\bk=\bff_{60}$;
$d=62$, then $\bk=\bff_{62}$;
$d=66$, then $\bk=\bff_{66}$;
$d=70$, then $\bk=\bff_{70}$;
$d=72$, then $\bk=\bff_{72}$;
$d=78$, then $\bk=\bff_{78}$;
$d=80$, then $\bk=\bff_{80}$;
$d=84$, then $\bk=\bff_{84}$;
$d=90$, then $\bk=\bff_{90}$;
$d=96$, then $\bk=\bff_{96}$;
$d=102$, then $\bk=\bff_{102}$;
$d=120$, then $\bk=\bff_{120}$;

\medskip

For $t=16^2$, we obtain that only the following $d\ge 6$ are possible,
and we get the corresponding estimates for $[\bk:\bq]$:

\noindent
$d=6$, then $[\bk:\bq]\le 8$;
$d=7$, then $[\bk:\bq]\le 144$;
$d=8$, then $[\bk:\bq]\le 18$;
$d=9$, then $[\bk:\bq]\le 21$;
$d=10$, then $[\bk:\bq]\le 10$;
$d=11$, then $[\bk:\bq]\le 30$;
$d=12$, then $[\bk:\bq]\le 10$;
$d=13$, then $[\bk:\bq]\le 24$;
$d=14$, then $[\bk:\bq]\le 9$;
$d=15$, then $[\bk:\bq]\le 8$;
$d=16$, then $[\bk:\bq]\le 12$;
$d=17$, then $[\bk:\bq]\le 24$;
$d=18$, then $[\bk:\bq]\le 12$;
$d=19$, then $[\bk:\bq]\le 18$;
$d=20$, then $[\bk:\bq]\le 12$;
$d=21$, then $[\bk:\bq]\le 12$;
$d=22$, then $[\bk:\bq]\le 10$;
$d=23$, then $\bk=\bff_{23}$;
$d=24$, then $[\bk:\bq]\le 12$;
$d=25$, then $\bk=\bff_{25}$;
$d=26$, then $[\bk:\bq]\le 12$;
$d=27$, then $\bk=\bff_{27}$;
$d=28$, then $[\bk:\bq]\le 12$;
$d=29$, then $\bk=\bff_{29}$;
$d=30$, then $[\bk:\bq]\le 12$;
$d=31$, then $\bk=\bff_{31}$;
$d=32$, then $\bk=\bff_{32}$;
$d=33$, then $\bk=\bff_{33}$;
$d=34$, then $\bk=\bff_{34}$;
$d=35$, then $\bk=\bff_{35}$;
$d=36$, then $[\bk:\bq]\le 12$;
$d=37$, then $\bk=\bff_{37}$;
$d=38$, then $\bk=\bff_{38}$;
$d=39$, then $\bk=\bff_{39}$;
$d=40$, then $\bk=\bff_{40}$;
$d=42$, then $[\bk:\bq]\le 12$;
$d=44$, then $\bk=\bff_{44}$;
$d=45$, then $\bk=\bff_{45}$;
$d=46$, then $\bk=\bff_{46}$;
$d=48$, then $\bk=\bff_{48}$;
$d=50$, then $\bk=\bff_{50}$;
$d=52$, then $\bk=\bff_{52}$;
$d=54$, then $\bk=\bff_{54}$;
$d=56$, then $\bk=\bff_{56}$;
$d=58$, then $\bk=\bff_{58}$;
$d=60$, then $\bk=\bff_{60}$;
$d=62$, then $\bk=\bff_{62}$;
$d=64$, then $\bk=\bff_{64}$;
$d=66$, then $\bk=\bff_{66}$;
$d=68$, then $\bk=\bff_{68}$;
$d=70$, then $\bk=\bff_{70}$;
$d=72$, then $\bk=\bff_{72}$;
$d=78$, then $\bk=\bff_{78}$;
$d=80$, then $\bk=\bff_{80}$;
$d=84$, then $\bk=\bff_{84}$;
$d=90$, then $\bk=\bff_{90}$;
$d=96$, then $\bk=\bff_{96}$;
$d=102$, then $\bk=\bff_{102}$;
$d=120$, then $\bk=\bff_{120}$;

To improve these results, let us apply a similar polynomial to the one we used
for $d=5$. For $d\ge 5$, we consider the integral polynomial
\begin{equation}
P(x)=x(x-1)(x+1-4\sin^2{\frac{\pi}{d}})\left(x-4\sin^2{\frac{\pi}{d}}\right)
\label{polynPd}
\end{equation}
over the ring of integers of $\bff_d$.
For $\sigma:\bff_d\to \br$ we denote by $\delta(\sigma)$ the maximum of
$|P^\sigma(x)|$ on the interval $[0,\sigma(4\sin^2{(\pi/d)}]$
(by calculating zeros of the derivative). Remark that
all $\varphi(d)/2$ embeddings $\sigma$ are defined by
$\sigma(4\sin^2{(\pi/d)})=4\sin^2{(k\pi/d)}$
where $1\le k\le d/2$ and $(k,d)=1$, and $\sigma^{(+)}$ corresponds to $k=1$.

Since evidently $P(\alpha)\not=0$, by Lemma \ref{lempolmain}, we get
\begin{equation}
1\le m\le \frac{\ln{P^{\sigma^{(+)}}(t^2)}-\ln{\delta(\sigma^{(+)})}}
{-\ln{\left(\prod_\sigma{\delta(\sigma)}\right)}}
\label{boundforPd}
\end{equation}
if $\prod_\sigma{\delta(\sigma)}<1$.
It usually gives a better upper bound for $m$ than \eqref{estfor3},
but it is harder to calculate.
Calculating \eqref{boundforPd} for $t=14^2$, $t=16^2$ and for
$d$ from the lists above, for some $d$ we get  a better upper bound for $m$.
We formulate the final result in the theorems below.

\begin{theorem}
For $d\ge 3$, V-arithmetic connected edge graphs
$\Gamma_2^{(3)}(14)$ (see Figure \ref{graphg3}) of the minimality $14$
(equivalently, totally real algebraic integers
$\alpha=u^2$ over $\bff_d=\bq(\cos^2(\pi/d))$ satisfying the conditions
\eqref{F32cond1} and \eqref{F32cond2} above for $t=14$)
are possible only for $d$ which are shown below. Moreover, we have
the shown below upper bounds for the degree $m=[\bk:\bff_d]$
of the ground field $\bk=\bq(\alpha)\supset \bff_d$:

\noindent
$d=3$, $m\le 44$; $d=4$, $m\le 16$; $d=5$, $m\le 10$;
$d=6$, $m\le 8$; $d=7$, $m\le 9$;
$d=8$, $m\le 8$; $d=9$, $m\le 5$; $d=10$, $m\le 5$;
$d=11$, $m\le 4$; $d=12$, $m\le 5$; $d=13$, $m\le 3$;
$d=14$, $m\le 3$; $d=15$, $m\le 2$; $d=16$, $m\le 3$;
$d=17$, $m\le 2$, $d=18$, $m\le 3$;
$d=20$, $m\le 2$, $d=21$, $m\le 2$; $d=22$, $m\le 2$;
$d=24$, $m\le 3$; $d=26$, $m\le 2$; $d=28$, $m\le 2$;
$d=30$, $m\le 3$; $d=36$, $m\le 2$; $d=42$, $m\le 2$;
for all remaining $d=19$, $23$, $25$, $27$, $29$, $31\ -\ 35$,\
$38\ -\ 40$,
$44$, $45$, $46$, $48$, $50$, $52$, $54$, $56$, $58$, $60$,
$66$, $70$,
$72$, $78$, $84$, $90$, $102$, $120$ we have
$m=1$ and $\bk=\bff_d$.

In particular, $[\bk:\bq]\le 44$ for $d=3$; $[\bk:\bq]\le 27$ for $d=7$;
$[\bk:\bq]\le 20$ for $d=5$, $11$; $[\bk:\bq]\le 18$ for $d=13$;
$[\bk:\bq]\le 16$ for all remaining $d$.
\label{thforF32and14}
\end{theorem}

\begin{theorem}
For $d\ge 3$, V-arithmetic connected edge graphs
$\Gamma_2^{(3)}(14)$ (see Figure \ref{graphg3}) of the minimality $16$
(equivalently, totally real algebraic integers
$\alpha=u^2$ over $\bff_d=\bq(\cos^2(\pi/d))$ satisfying the conditions
\eqref{F32cond1} and \eqref{F32cond2} above for $t=16$)
are possible only for $d$ which are shown below. Moreover, we have
the shown below upper bounds for the degree $m=[\bk:\bff_d]$
of the ground field $\bk=\bq(\alpha)\supset \bff_d$:

Then only the following $d$ are possible, and we have
the following upper bounds for $m=[\bk:\bff_d]$:

\noindent
$d=3$, $m\le 47$; $d=4$, $m\le 16$; $d=5$, $m\le 10$;
$d=6$, $m\le 8$; $d=7$, $m\le 10$;
$d=8$, $m\le 9$; $d=9$, $m\le 5$; $d=10$, $m\le 5$;
$d=11$, $m\le 4$; $d=12$, $m\le 5$; $d=13$, $m\le 3$;
$d=14$, $m\le 3$; $d=15$, $m\le 2$; $d=16$, $m\le 3$;
$d=17$, $m\le 2$, $d=18$, $m\le 4$; $d=19$, $m\le 2$;
$d=20$, $m\le 3$, $d=21$, $m\le 2$; $d=22$, $m\le 2$;
$d=24$, $m\le 3$; $d=26$, $m\le 2$; $d=28$, $m\le 2$;
$d=30$, $m\le 3$; $d=36$, $m\le 2$; $d=42$, $m\le 2$;
for all remaining $d=23$, $25$, $27$, $29$, $31\ -\ 35$,\
$38\ -\ 40$,
$44$, $45$, $46$, $48$, $50$, $52$, $54$, $56$, $58$, $60$, $62$,
$64$, $66$, $70$,
$72$, $78$, $80$, $84$, $90$, $96$, $102$, $120$ we have
$m=1$ and $\bk=\bff_d$.

In particular, $[\bk:\bq]\le 47$ for $d=3$; $[\bk:\bq]\le 30$ for $d=7$;
$[\bk:\bq]\le 20$ for $d=5$, $11$; $[\bk:\bq]\le 18$ for $d=8$, $13$, $19$;
$[\bk:\bq]\le 16$ for all remaining $d$.

\label{thforF32and16}
\end{theorem}

\medskip

By Theorem \ref{thforF32and14}, we obtain

\begin{theorem} The set of fields $\F \Gamma_2^{(3)}(14)$ is finite and
their degree over $\bq$ has the upper bound $44$. See more exact information
about these fields depending on the parameter $d$ (see Figure \ref{graphg3})
in Theorem \ref{thforF32and14}.
\label{thfieldsGamma2(3)}
\end{theorem}

\subsubsection{Ground fields of some $\Gamma_1^{(3)}(t)$}
\label{subsubsec:Gamma1(3)}
Here we consider ground fields $\F \Gamma_1^{(3)}(t)$ of V-arithmetic
graphs $\Gamma_1^{(3)}(t)$ for parameters $3\le  s\le k\le 5$.
See Figure \ref{graphg3}.

The corresponding Gram matrix is
\begin{equation}
\left(\begin{array}{ccc}
-2     & u                    &   2\cos{\frac{\pi}{s}}                   \\
 u     & -2                   & 2\cos{\frac{\pi}{k}}  \\
 2\cos{\frac{\pi}{s}}     &  2\cos{\frac{\pi}{k}}  & -2
\end{array}\right)
\label{gramm1(3)}
\end{equation}
where $s,k\ge 3$ are integers, and $u$ is a totally real algebraic integer.
The ground field is
$$
\bk=\bq\left(u^2,\cos^2{(\pi/s)}, \cos^2{(\pi/k)},
u\cos{(\pi/s)}\cos{(\pi/k)}\right).
$$
The determinant $d(u)$ of the Gram matrix is given by the equality
$$
\frac{d(u)}{2}=u^2+4\cos{\frac{\pi}{s}}\cos{\frac{\pi}{k}}u+
4\cos^2{\frac{\pi}{s}}+4\cos^2{\frac{\pi}{k}}-4.
$$
It follows that $\Gamma_1^{(3)}$ is V-arithmetic of the minimality $t>2$,
if and only if we have
$$
0<-\widetilde{\sigma}(2\cos{\frac{\pi}{s}}\cos{\frac{\pi}{k}})-
\sqrt{\sigma(4\sin^2{\frac{\pi}{s}}\sin^2{\frac{\pi}{k}})}<
\widetilde{\sigma}(u)<
$$
\begin{equation}
<-\widetilde{\sigma}(2\cos{\frac{\pi}{s}}\cos{\frac{\pi}{k}})+
\sqrt{\sigma(4\sin^2{\frac{\pi}{s}}\sin^2{\frac{\pi}{k}})}
\label{F31cond1}
\end{equation}
for all $\sigma:\bk\to \br$ which are different from $\sigma^{(+)}$, and
\begin{equation}
4<\sigma^{(+)}(u^2)<t^2.
\label{F31cond2}
\end{equation}
where $\widetilde{\sigma}$ extends $\sigma:\bk=\bq(u^2)\to \br$
to $\bq(u)=\bk(\cos{(\pi/s)}\cos{(\pi/k)})$.

\medskip

{\bf Case $s=k=3$.} In this case, $\bk=\bq(u)$ and
$$
-2<\sigma(u)<1,\ \ \ \ 2<\sigma^{(+)}(u)<t.
$$
We take the polynomial
\begin{equation}
P(x)=(x+2)^3(x+1)^4x^4(x-1)^3(x^2+x-1)^2
\label{polGamma31for3,3}
\end{equation}
which is similar to \eqref{pol3}, and we similarly get
$$
[\bk:\bq]\le 1+\frac{\ln(P(t))}{(-\ln(\delta))}
$$
where $\delta= 0.0905969664$. For $t=14$, we get
$[\bk:\bq]\le 23$. (Theorems \ref{th121} and \ref{th122} give only
the upper bound $57$.)

\medskip

{\bf Case $s=3$, $k=4$.} Then $\bq(\sqrt{2})\subset \bk(u)$ and
$\bk=\bq(u^2)$. We have $\sqrt{2}\in \bk$ if and only if $u\in \bk$.
In this case,
$$
d(u)/2=u^2+\sqrt{2}u-1.
$$
We have
$$
\widetilde{\sigma}(-\frac{1}{\sqrt{2}})-\sqrt{\frac{3}{2}}<
\widetilde{\sigma}(u)<
\widetilde{\sigma}(-\frac{1}{\sqrt{2}})+\sqrt{\frac{3}{2}},
$$
$$
4<\sigma^{(+)}(u^2)<t^2.
$$
We take the polynomial
\begin{equation}
P(x)=x^4(x+\sqrt{2})^4(x^2+\sqrt{2}x-1)
\label{polGamma31for3,4}
\end{equation}
over the ring of integers of $\bq(\sqrt{2})$. The maximum of
$P(x)$ on the interval $[-1/\sqrt{2}-\sqrt{3/2},-1/\sqrt{2}+\sqrt{3/2}]$
is equal to $\delta_1=0.09375...$.
The maximum of the conjugate polynomial
$P^g(x)=x^4(x-\sqrt{2})^4(x^2-\sqrt{2}x-1)$ on the interval
$[1/\sqrt{2}-\sqrt{3/2},1/\sqrt{2}+\sqrt{3/2}]$ is equal to
$\delta_2=0.09375...$. By Lemma \ref{lempolmain}, we get
$$
[\bk:\bq(\sqrt{2})]<\frac{\ln{P(t)}-\ln{\delta_1}}{-\ln(\delta_1\delta_2)}
$$
if $\sqrt{2}\in \bk$, and
$$
[\bk(u):\bq(\sqrt{2})]<\frac{\ln{((P(t)P^g(-t))}-
\ln{\delta_1\delta_2}}{-log(\delta_1\delta_2)}
$$
if $\sqrt{2}$ does not belong to $\bk$.

For $t=14$, in both cases we get $[\bk:\bq]\le 12$.

\medskip

{\bf Case $s=k=4$.} Then $u\in \bk$ and $\bk=\bq(u)$. We have
$$
-2<\sigma(u)<0,\ \ \ \ 2<\sigma^{(+)}(u)<14.
$$
We take
\begin{equation}
P(x)=(x+2)(x+1)^2x
\label{polGamma31for4,4}
\end{equation}
which is similar to the polynomial \eqref{pol4}. 
By Lemma \ref{lempolmain}, 
$$
[\bk:\bq]\le 1+\frac{\ln{P(t)}}{-\ln{(1/4)}}.
$$
For $t=14$, we obtain $[\bk:\bq]\le 8$.


\medskip

{\bf Case $s=3$, $k=5$.} Then $u\in \bk$ and
$\bff_5=\bq(\cos^2{(\pi/5)})\subset \bk=\bq(u)$. We have
$$
\frac{d(u)}{2}=u^2+2\cos{\frac{\pi}{5}}u+4\cos^2{\frac{\pi}{5}}-3\ .
$$
We have
\begin{equation}
-\sigma(\cos{\frac{\pi}{5}})-\sqrt{3}\,\sigma(\sin{\frac{\pi}{5}})<\sigma(u)
<-\sigma(\cos{\frac{\pi}{5}})+\sqrt{3}\,\sigma(\sin{\frac{\pi}{5}}),\ \
2<\sigma^{(+)}(u)<14.
\label{Gamma32int3,5}
\end{equation}
We take the polynomial

$$
P(x)=
$$
\begin{equation}
x(x+1)(x+3-4\sin^2{(\pi/5)})(x+2-4\sin^2{(\pi/5)})
(x^2+2\cos{\frac{\pi}{5}}x+4\cos^2{\frac{\pi}{5}}-3)
\label{polGamma31for3,5}
\end{equation}
over the ring of integers of $\bff_5$.
The maxima $\delta(\sigma)$ of $|P^\sigma(x)|$ 
on the corresponding intervals 
\eqref{Gamma32int3,5} are equal to $\delta_1=0.0690098...$ if
$\sigma(\sin{(\pi/5)})=\sin{(\pi/5)}$, and $\delta_2=1.2383316...$ if
 $\sigma(\sin{(\pi/5)})=\sin{(3\pi/5)}$. By Lemma \ref{lempolmain},
$$
[\bk:\bff_5]<\frac{\ln{P(t)}-\ln{\delta_1}}{-\ln(\delta_1\delta_2)}.
$$
For $t=14$, we get $[\bk:\bff_5]\le 7$ and $[\bk:\bq]\le 14$.

\medskip

\medskip

{\bf Case $s=4$, $k=5$.} (This case is the most difficult.)\newline
Then $\bff_5=\bq(\cos^2(\pi/5))\subset
\bk =\bq(u^2)$, and
$u\in \bk$ if and only if $\sqrt{2}\in \bk$. We have
$$
\frac{d(u)}{2}=u^2+2\sqrt{2}\cos{\frac{\pi}{5}}\,u+4\cos^2{\frac{\pi}{5}}-2=
$$
$$
u^2+\left(\frac{\sqrt{2}}{2}+\frac{\sqrt{10}}{2}\right)u+
\left(-\frac{1}{2}+\frac{\sqrt{5}}{2}\right),
$$
and
$$
0<a(\widetilde{\sigma})=-\widetilde{\sigma}(\sqrt{2}\cos{\frac{\pi}{5}})-
\sqrt{\sigma(2\sin^2{\frac{\pi}{5}})}<
\widetilde{\sigma}(u)<
$$
\begin{equation}
<-\widetilde{\sigma}(\sqrt{2}\cos{\frac{\pi}{5}})+
\sqrt{\sigma(2\sin^2{\frac{\pi}{5}})}=b(\widetilde{\sigma})
\label{F32cond1for4,5}
\end{equation}
for all $\sigma:\bk\to \br$ which are different from $\sigma^{(+)}$, and
\begin{equation}
4<\sigma^{(+)}(u^2)<t^2.
\label{F32cond2for4,5}
\end{equation}
We take the polynomial $P(x)$ of the degree $8$,
$$
P(x)=\left(x^2+\left(\frac{\sqrt{2}}{2}+\frac{\sqrt{10}}{2}\right)x-
\frac{1}{2}+\frac{\sqrt{5}}{2}\right)\times
$$
$$
\left(4x^3+(3\sqrt{2}+3\sqrt{10})x^2+(6+4\sqrt{5})x+\frac{5\sqrt{2}}{2}+
\frac{\sqrt{10}}{2}\right)\times
$$
\begin{equation}
\left(4x^3+(3\sqrt{2}+3\sqrt{10})x^2+(7+3\sqrt{5})x+\frac{3\sqrt{2}}{2}+
\frac{\sqrt{10}}{2}\right)
\label{polGamma31for4,5}
\end{equation}
over the ring of integers of $\bff_{4,5}=\bq(\sqrt{2},\cos^2(\pi/5))=
\bq(\sqrt{2},\sqrt{5})$.

All roots of the polynomials $P^{\widetilde{\sigma}}(x)$ belong to the
corresponding intervals $[a(\widetilde{\sigma}),\,b(\widetilde{\sigma})]$.
The maxima $\delta(\widetilde\sigma)$ of $|P^{\widetilde{\sigma}}(x)|$ on the
corresponding intervals \eqref{F32cond1for4,5} are as follows.
The maxima of $|P^{\widetilde{\sigma}}(x)|$ on the interval
$$
[a(\widetilde{\sigma}),\ 
b(\widetilde{\sigma})]=[-1.97537668..,\ -0.31286893...]
$$
is equal to $\delta_{1}=0.045593135...$ 
if $\widetilde{\sigma}(\sqrt{2})=\sqrt{2}$
and $\widetilde{\sigma}(\sqrt{5})=\sqrt{5}$;
on the interval
$$
[a(\widetilde{\sigma}),\ 
b(\widetilde{\sigma})]=[0.31286893...,\-1.97537668...]
$$
it is equal to $\delta_{2}=0.045593135...$
if $\widetilde{\sigma}(\sqrt{2})=-\sqrt{2}$ and
$\widetilde{\sigma}(\sqrt{5})=\sqrt{5}$;
on the interval
$$
[a(\widetilde{\sigma}),\ b(\widetilde{\sigma})]=
[-0.90798099...,\ 1.78201304...]
$$
it is equal to $\delta_{3}=2.14190686...$
if $\widetilde{\sigma}(\sqrt{2})=\sqrt{2}$ and
$\widetilde{\sigma}(\sqrt{5})=-\sqrt{5}$;
on the interval
$$
[a(\widetilde{\sigma}),\ b(\widetilde{\sigma})]=
[-1.78201304...,\ 0.90798099...]
$$
it is equal to $\delta_{4}=2.14190686...$
if $\widetilde{\sigma}(\sqrt{2})=-\sqrt{2}$ and
$\widetilde{\sigma}(\sqrt{5})=-\sqrt{5}$.

By Lemma \ref{lempolmain},
$$
[\bk:\bff_{4,5}]\le \frac{\ln{P(t)}-\ln{\delta_{1}}}
{-\ln{(\delta_{1}\delta_{2}\delta_{3}\delta_{4})}}
$$
if $\sqrt{2}\in \bk$,
and
$$
[\bk(u):\bff_{4,5}]\le \frac{\ln{(P(t)P^\tau(-t))}-\ln{(\delta_{1}\delta_{2})}}
{-\ln{(\delta_{1}\delta_{2}\delta_{3}\delta_{4})}}
$$
if $\sqrt{2}\notin \bk$ where $\tau(\sqrt{2})=-\sqrt{2}$ and
$\tau(\sqrt{5})=\sqrt{5}$. For $t=14$, in the first case, we get
$[\bk:\bff_{4,5}]\le 5$ and $[\bk:\bq]\le 20$; for the second case,
we get $[\bk(u):\bff_{4,5}]\le 11$ and $[\bk(u):\bq]\le 44$.
Then $[\bk:\bq]\le 22$. Thus, in both cases, we get $[\bk:\bq]\le 22$.
(Theorems \ref{th121} and \ref{th122} give $[\bk:\bq]\le 72$ only.)

\medskip

{\bf Case $s=k=5$.} Then $\bff_5=\bq(\cos^2{(\pi/5)})\subset \bk=\bq(u)$,
and we have
\begin{equation}
-2<\sigma(u)<\sigma(4\sin^2{\frac{\pi}{5}})-2,\ \ 2<\sigma^{(+)}(u)<14.
\label{Gamma32int5,5}
\end{equation}
We take the polynomial
\begin{equation}
P(x)=(x+2)(x+1)(x+3-4\sin^2{\frac{\pi}{5}})(x+2-4\sin^2{\frac{\pi}{5}})
\label{polGamma31for5,5}
\end{equation}
which is similar to \eqref{pol5}. The maxima $\delta(\sigma)$
of $|P^\sigma(x)|$ on the 
corresponding intervals \eqref{Gamma32int5,5}
are equal to $\delta_1=0.0455931...$ 
if $\sigma(\sin^2{(\pi/5)})=\sin^2{(\pi/5)}$;
$\delta_2=2.1419068...$ if $\sigma(\sin^2{(\pi/5)})=\sin^2{(2\pi/5)}$.

By Lemma \ref{lempolmain},
$$
[\bk:\bff_5]\le \frac{\ln{P(t)}-\ln{\delta_1}}{-\ln{(\delta_1\delta_2)}}.
$$
For $t=14$, we get $[\bk:\bff_5]\le 6$ and $[\bk:\bq]\le 12$.

Thus, finally we get the following statement.

\begin{theorem} We have the following upper bounds for the degrees
of ground fields $\bk\in \F \Gamma_1^{(3)}(14)$ that is for
V-arithmetic edge graphs
  $\Gamma_1^{(3)}$ of the minimality $14$ (i.e., $2<u<14$)
with parameters $3\le s\le k\le 5$ (see Figure \ref{graphg3}):

\noindent
$[\bk:\bq]\le 23$ for $s=k=3$;

\noindent
$[\bk:\bq]\le 12$ for $s=3$, $k=4$;

\noindent
$[\bk:\bq]\le 8$ for $s=k=4$;

\noindent
$\bff_5=\bq(\sqrt{5})\subset \bk$ and $[\bk:\bff_5]\le 7$ for $s=3$, $k=5$;

\noindent
for $s=4$, $k=5$ we have
$\bff_{4,5}=\bq(\sqrt{2},\sqrt{5})\subset \bk$ and $[\bk:\bff_{4,5}]\le 5$
if $\sqrt{2}\in \bk$, and $\bff_5=\bq(\sqrt{5})\subset \bk$ and
$[\bk:\bff_5]\le 11$ if $\sqrt{2}\notin \bk$;

\noindent
$\bff_5=\bq({\sqrt{5}})\subset \bk$ and $[\bk:\bff_5]\le 6$ for $s=k=5$.

In particular, $[\bk:\bq]\le 23$ for all fields
$\bk\in \F \Gamma_1^{(3)}(14)$ with parameters $3\le s\le k\le 5$.

\label{thforF31and14}
\end{theorem}


\subsection{Fields from $\F \Gamma_1^{(4)}(14)$}\label{subsec:gam(4)1}

\begin{figure}
\begin{center}
\includegraphics[width=6cm]{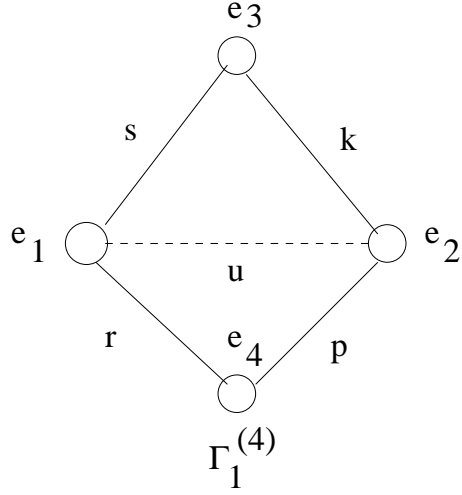}
\end{center}
\caption{The graph $\Gamma_1^{(4)}$}
\label{graphg41}
\end{figure}

For $\Gamma_1^{(4)}(14)$ (see Figure \ref{graphg41}) we assume that integers
$s,k,r,p\ge 3$. Subgraphs
$\Gamma_1^{(4)}-\{e_1\}$ and $\Gamma_1^{(4)}-\{e_2\}$ must be Coxeter
graphs. It follows that
we must consider only (up to obvious symmetries) the following cases:
either $s=k=3$ and $5\ge r\ge p\ge 3$, or $s=p=3$ and $5\ge r\ge k\ge 4$;
the totally real algebraic integer $u$ satisfies the inequality
$2<u<14$.

Since the graph $\Gamma_1^{(4)}(14)$ with the parameters $s,\,k,\,r,\,p$
contains two subraphs $\Gamma_1^{(3)}(14)$ with the parameters
$s,k$ and $r,p$, we obtain that any field $\bk\in \F \Gamma_1^{(4)}(14)$
with the parameters $s,k,r,p$ belongs to
$\F \Gamma_1^{(3)}(14)$ with the parameters $s,k$ and to
$\F \Gamma_1^{(3)}(14)$ with the parameters $p,r$. Thus, the degree of
fields $\bk\in \F \Gamma_1^{(4)}(14)$
with the parameters $s,k,r,p$ is less than minimum of degrees of fields
from  $\F \Gamma_1^{(3)}(14)$ with the parameters $s,k$ and fields from
$\F \Gamma_1^{(3)}(14)$ with the parameters $p,r$. Thus, applying
Theorem \ref{thforF31and14} we obtain an ``easy'' result.

\begin{proposition}
 We have the following upper bounds for the degrees
of ground fields $\bk\in \F \Gamma_1^{(4)}(14)$ that is for
V-arithmetic edge graphs
  $\Gamma_1^{(4)}$ of the minimality $14$ (i.e., $2<u<14$) depending
on their parameters $s,k,p,r$ (see Figure \ref{graphg41}):

\noindent
$[\bk:\bq]\le 23$ for $s=k=p=r=3$;

\noindent
$[\bk:\bq]\le 12$ for $s=k=3$, $r=4$, $p=3$;

\noindent
$\bff_5=\bq(\sqrt{5})\subset \bk$ and
$[\bk:\bff_5]\le 7$ for $s=k=3$, $r=5$, $p=3$;

\noindent
$[\bk:\bq]\le 8$ for $s=k=3$, $r=p=4$;

\noindent
$\bff_{4,5}=\bq(\sqrt{2},\sqrt{5})\subset \bk$ and
$[\bk:\bff_{4,5}]\le 5$ if $s=k=3$, $r=5$, $p=4$;

\noindent
$\bff_5=\bq({\sqrt{5}})\subset \bk$ and $[\bk:\bff_5]\le 6$ for
for $s=k=3$, $r=p=5$.

\noindent
$[\bk:\bq]\le 12$ for $s=p=3$, $r=k=4$;

\noindent
$\bff_{4,5}=\bq(\sqrt{2},\sqrt{5})\subset \bk$ and
$[\bk:\bff_{4,5}]\le 3$ for $s=p=3$, $r=5$, $k=4$;

\noindent
$\bff_5=\bq({\sqrt{5}})\subset \bk$ and $[\bk:\bff_5]\le 6$ for
$s=p=3$, $r=k=5$.

In particular, $[\bk:\bq]\le 23$ for all fields
$\bk\in \F \Gamma_1^{(4)}(14)$.

\label{propforF41and14}
\end{proposition}

By considering the graphs $\Gamma_1^{(4)}$ directly, one can 
significantly improve Proposition \ref{propforF41and14}. 
We hope to consider this in further variants of the paper and 
further publications.

\subsection{Fields from $\F \Gamma_2^{(4)}(14)$}\label{subsec:gam(4)2}
For $\Gamma_2^{(4)}(14)$ (see Figure \ref{graphg42}),
$s,k,p \ge 3$ are natural numbers and $2<u<14$ is
a totally real algebraic integer. Moreover, we have only the following
possibilities: $3\le s\le k\le 5$, $p=3$; $s=k=3$, $p=4,\,5$.

Since the graph $\Gamma_2^{(4)}(14)$ with the parameters $s,\,k,\,p$
contains the subraph $\Gamma_1^{(3)}(14)$ with the parameters
$s,k$, any field $\bk\in \F \Gamma_2^{(4)}(14)$
with the parameters $s,k,p$ belongs to
$\F \Gamma_1^{(3)}(14)$ with the parameters $s,k$. 
Thus, the degree of
fields $\bk\in \F \Gamma_2^{(4)}(14)$
with the parameters $s,k,p$ is less than the degrees of fields
from  $\F \Gamma_1^{(3)}(14)$ with the parameters $s,k$. Thus, applying
Theorem \ref{thforF31and14} we obtain an ``easy'' result.

\begin{proposition}
 We have the following upper bounds for the degrees
of ground fields $\bk\in \F \Gamma_2^{(4)}(14)$ that is for
V-arithmetic edge graphs
  $\Gamma_2^{(4)}$ of the minimality $14$ (i.e., $2<u<14$) depending
on their parameters $s,k,p$ (see Figure \ref{graphg42}):

\noindent
$[\bk:\bq]\le 23$ for $s=k=3$, $p=3,4,5$;

\noindent
$[\bk:\bq]\le 12$ for $s=3$, $k=4$, $p=3$;

\noindent
$[\bk:\bq]\le 8$ for $s=k=4$, $p=3$;

\noindent
$\bff_5=\bq(\sqrt{5})\subset \bk$ and $[\bk:\bff_5]\le 7$ 
for $s=3$, $k=5$, $p=3$;

\noindent
for $s=4$, $k=5$, $p=3$ we have
$\bff_{4,5}=\bq(\sqrt{2},\sqrt{5})\subset \bk$ and $[\bk:\bff_{4,5}]\le 5$
if $\sqrt{2}\in \bk$, and $\bff_5=\bq(\sqrt{5})\subset \bk$ and
$[\bk:\bff_5]\le 11$ if $\sqrt{2}\notin \bk$;

\noindent
$\bff_5=\bq({\sqrt{5}})\subset \bk$ and $[\bk:\bff_5]\le 6$ for $s=k=5$, 
$p=3$. 

In particular, $[\bk:\bq]\le 23$ for all fields
$\bk\in \F \Gamma_2^{(4)}(14)$.

\label{propforF42and14}
\end{proposition}

By considering the graphs $\Gamma_2^{(4)}$ directly, one can 
significantly improve Proposition \ref{propforF42and14}. 
We hope to consider this in further variants of the paper and 
further publications.

\begin{figure}
\begin{center}
\includegraphics[width=6cm]{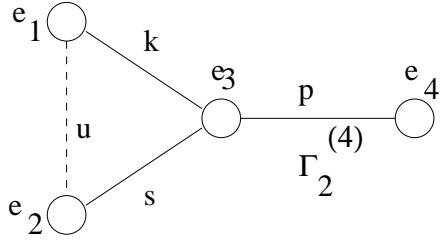}
\end{center}
\caption{The graph $\Gamma_2^{(4)}$}
\label{graphg42}
\end{figure}

\subsection{Fields from $\F \Gamma_3^{(4)}(14)$}\label{subsec:gam(3)4}
For $\Gamma_3^{(4)}(14)$ (see Figure \ref{graphg43}),
$s\ge 2$, $k,r \ge 3$ are natural numbers, and $2<u<14$ is
a totally real algebraic integer. Moreover, we have only the following
possibilities: $s=2$, $k=3$, $r=3,\,4,\,5$;
$s=2$, $k=4,\,5$, $r=3$; $3\le s\le k\le 5$, $r=3$; $s=k=3$, $r=4,\,5$.

The ground field $\bk=\bq(u^2)$ contains cyclic products
$$
\cos^2{\frac{\pi}{s}},\ \cos^2{\frac{\pi}{k}},\ \cos^2{\frac{\pi}{r}},\
u^2,\ u\,\cos{\frac{\pi}{s}}\cos{\frac{\pi}{k}}\cos{\frac{\pi}{r}}\,.
$$ 
The determinant $d(u)$  of the Gram matrix is determined by the equality
\begin{equation}
-\frac{d(u)}{4}=
\sin^2{\frac{\pi}{r}}\,u^2+2\cos{\frac{\pi}{s}} \cos{\frac{\pi}{k}}
\cos{\frac{\pi}{r}}\,u+4 \cos^2{\frac{\pi}{r}}-
4 \sin^2{\frac{\pi}{s}}\sin^2{\frac{\pi}{k}}.
\label{detforF43}
\end{equation}
Let
$$
D=4\cos^2{\frac{\pi}{s}} \cos^2{\frac{\pi}{k}} \cos^2{\frac{\pi}{r}}+
16 \sin^2{\frac{\pi}{s}}\sin^2{\frac{\pi}{k}}\sin^2{\frac{\pi}{r}}-
16 \sin^2{\frac{\pi}{r}}\cos^2{\frac{\pi}{r}}
$$
be the discriminant of this quadratic polynomial of the variable $u$.
The graph $\Gamma^{(4)}_3(14)$ is V-arithmetic if and only if for
$\tau:\bk(u)\to \br$ which is different from $\sigma^{(+)}$ on $\bk$,
one has
$$
\frac{-2\widetilde\tau\left(\cos{\frac{\pi}{s}} \cos{\frac{\pi}{k}}
\cos{\frac{\pi}{r}}\right)-\sqrt{\tau(D)}}{2\tau(\sin^2{\frac{\pi}{r}})}\
<\tau(u)<\
$$
\begin{equation}
< \frac{-2\widetilde\tau\left(\cos{\frac{\pi}{s}} \cos{\frac{\pi}{k}}
\cos{\frac{\pi}{r}}\right)+\sqrt{\tau(D)}}{2\tau(\sin^2{\frac{\pi}{r}})}
\label{infor43and14}
\end{equation}
where $\widetilde{\tau}$ extends $\tau$.

\begin{figure}
\begin{center}
\includegraphics[width=6cm]{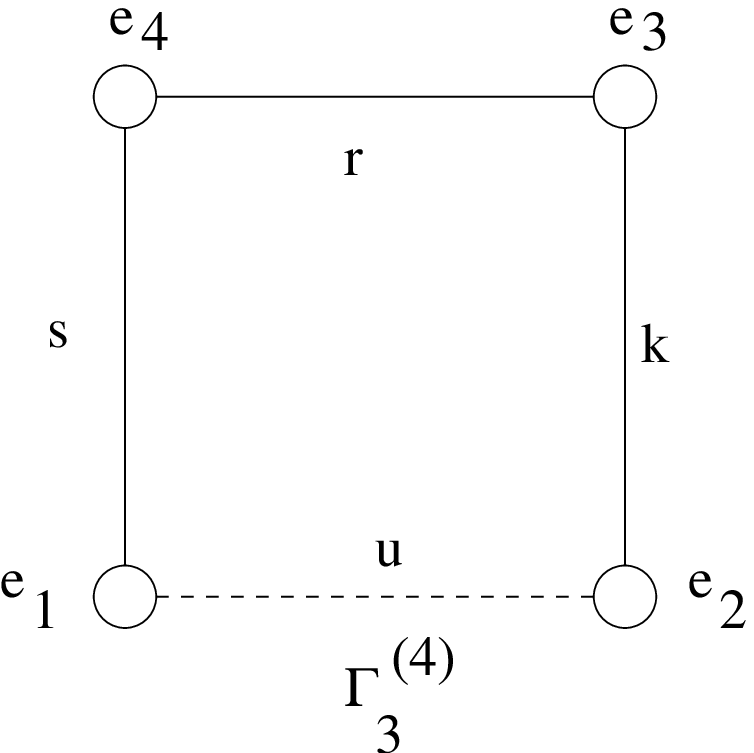}
\end{center}
\caption{The graph $\Gamma_3^{(4)}$}
\label{graphg43}
\end{figure}

Since the graph $\Gamma_3^{(4)}(14)$ with the parameters $s,\,k,\,r$ 
contains subraphs \linebreak $\Gamma_2^{(3)}(14)$ with the parameters
$s$, if $s\ge 3$,  and $3\le k\le 5$, 
we obtain that any field $\bk\in \F \Gamma_3^{(4)}(14)$
with the parameters $s,k,r$ belongs to
$\F \Gamma_2^{(3)}(14)$ with the parameter $s$, if $s\ge 3$,  
and to $\F \Gamma_2^{(3)}(14)$ with the parameter $3\le k\le 5$. 
Thus, the degree of fields $\bk\in \F \Gamma_3^{(4)}(14)$
with the parameters $s,k,r$ is less than minimum of degrees of fields
from  $\F \Gamma_2^{(3)}(14)$ with the parameter $s$, if $s\ge 3$,  
and fields from $\F \Gamma_2^{(3)}(14)$ with the parameters $3\le k\le 5$. 
Thus, applying Theorem \ref{thforF32and14} we obtain an ``easy'' result.

\begin{proposition}
 We have the following upper bounds for the degrees
of ground fields $\bk\in \F \Gamma_3^{(4)}(14)$ that is for
V-arithmetic edge graphs
  $\Gamma_3^{(4)}$ of the minimality $14$ (i.e., $2<u<14$) depending
on their parameters $s,k,r$ (see Figure \ref{graphg43}):

\noindent
$[\bk:\bq]\le 44$ for $s=2$, $k=3$, $r=3,4,5$;

\noindent
$[\bk:\bq]\le 16$ for $s=2$, $k=4$, $r=3$;

\noindent
$\bff_5=\bq(\sqrt{5})\subset \bk$ and 
$[\bk:\bff_5]\le 10$ for $s=2$, $k=5$, $r=3$; 

\noindent
$[\bk:\bq]\le 44$ for $s=k=3$, $r=3,4,5$;

\noindent
$[\bk:\bq]\le 16$ for $s=3$, $k=4$, $r=3$;

\noindent
$\bff_5=\bq(\sqrt{5})\subset \bk$ and 
$[\bk:\bff_5]\le 10$ for $s=3$, $k=5$, $r=3$;

\noindent
$[\bk:\bq]\le 16$ for $s=k=4$, $r=3$;

\noindent
$\bff_5=\bq(\sqrt{5})\subset \bk$ and 
$[\bk:\bff_5]\le 8$ for $s=4$, $k=5$, $r=3$;

\noindent
$\bff_5=\bq(\sqrt{5})\subset \bk$ and 
$[\bk:\bff_5]\le 10$ for $s=k=5$, $r=3$;

\label{propforF43and14}
\end{proposition}

Let us improve the poor bounds for $s=2$, $k=3$ and $s=k=3$.

For $s=2$, we have for $\alpha=u^2$ and $\sigma:\bk\to \br$ 
$$
0<\sigma(\alpha)<\sigma\left(D/(4\sin^4{(\pi/r)})\right)
$$
if $\sigma\not=\sigma^{(+)}$. 

\medskip 

For $s=2$, $k=r=3$, 
$$
0<\sigma (\alpha)< \frac{8}{3},\ \ \ 4<\sigma^{(+)}(\alpha)<14^2.
$$
We take the polynomial 
\begin{equation}
P(x)=x(x-1)^2(x-2)^2(x^2-3x+1)^2.
\label{pol41for2,3,3}
\end{equation}
The maximum of $|P(x)|$ on $[0,8/3]$ is $\delta=0.148148...$. By Lemma 
\ref{lempolmain}, 
$$
[\bk:\bq]\le 1+\frac{\ln{|P(14^2)|}}{-\ln{\delta}}\le 25.9. 
$$
Thus,  $[\bk:\bq]\le 25$.

\medskip

For $s=2$, $k=3$, $r=4$, 
$$
0<\sigma(\alpha)<2,\ \ \ 4<\sigma^{(+)}(\alpha)<14^2.
$$
We take the polynomial
\begin{equation}
P(x)=x(x-1)^2(x-2).
\label{pol14for2,3,4}
\end{equation}
The maximum of $|P(x)|$ on the interval $[0,2]$ is equal to
$\delta=1/4$. By Lemma \ref{lempolmain}, 
$$
[\bk:\bq]\le 1+\frac{\ln{P(14^2)}}{-\ln{\delta}}<16.3. 
$$
Thus, $[\bk:\bq]\le 16$.

\medskip 

For  $s=2$, $k=3$, $r=5$, 
$$
0<\sigma(\alpha)<4-(16/5)\sigma(\sin^2{(2\pi/5)}),\ \ \
4<\sigma^{(+)}(\alpha)<14^2.
$$
We take the polynomial
\begin{equation}
P(x)=x(x-1)(x+1-4\sin^2{(\pi/5)})
\label{pol14for2,3,5}
\end{equation}
over the ring of integers of the field
$\bff_5=\bq(\cos^2{(\pi/5)})$. The maximum
of $|P(x)|$ on $[0,4-(16/5)\sin^2{(2\pi/5)}]$ is equal to 
$\delta_1=0.0844582...$. The maximum of $|P^g(x)|$ on the interval 
$[0,4-(16/5)\sin^2{(\pi/5)}]$ where
$g:\bff_5\to \br$ is not identity, is equal to 
$\delta_2=1.51554175...$.
By Lemma \ref{lempolmain}, 
$$
[\bk:\bff_5]\le 
\frac{\ln{P(14^2)}-\ln {\delta_1}}{-\ln{(\delta_1\delta_2)}}\le 8.91.
$$
Thus, $[\bk:\bff_5]\le 8$. 

\medskip

Let $s=k=r=3$. Then 
$$
0<\sigma(u^2)<x_1^2=2.15612498...,\ \ \ 4<\sigma^{(+)}(u^2)<14^2
$$
where $x_1$ is the left root of the polynomial $3x^2+x-5$ 
defined by \eqref{detforF43}. 
We take the polynomial 
\begin{equation}
P(x)=x(x-1)(x-2). 
\label{pol14for3,3,3}
\end{equation}  
The maximum of $|P(x)|$ on the interval $[0,2.15612498...]$ 
is equal to \linebreak 
$\delta=0.3849001...$. By Lemma \ref{lempolmain},
$$
[\bk:\bq]\le 
1+\frac{\ln{P(14^2)}}{-\ln{\delta}}\le 17.2.
$$
Thus, $[\bk:\bq]\le 17$. 

\medskip

Let $s=k=3$, $r=4$. Then 
$$
0<\sigma(u^2)<x_1^2=1.30901699...,\ \ \ 4<\sigma^{(+)}(u^2)<14^2,
$$
where $x_1$ is the left root of the polynomial $2x^2+\sqrt{2}x-1$ 
defined by \eqref{detforF43}.
We take the polynomial
\begin{equation}
P(x)=x(x-1)^2. 
\label{pol14for3,3,4}
\end{equation}  
The maximum of $|P(x)|$ on the interval $[0,1.30901699...]$ 
is equal to \linebreak 
$\delta=0.14814...$. By Lemma \ref{lempolmain},
$$
[\bk:\bq]\le 
1+\frac{\ln{P(14^2)}}{-\ln{\delta}}\le 8.3.
$$
Thus, $[\bk:\bq]\le 8$. 

\medskip 

Let $s=k=3$, $r=5$. Then $\bff_5=\bq(\cos^2(\pi/5))\subset \bk=\bq(u)$. 
In this case, $\sigma(D)<0$ if $\sigma(\sin^2(\pi/5))=\sin^2(\pi/5)$.  
It follows that $\bk=\bff_5$. 

Thus, finally, we get the following result. 

\begin{theorem}
 We have the following upper bounds for the degrees
of ground fields $\bk\in \F \Gamma_3^{(4)}(14)$ that is for
V-arithmetic edge graphs
  $\Gamma_3^{(4)}$ of the minimality $14$ (i.e., $2<u<14$) depending
on their parameters $s,k,r$ (see Figure \ref{graphg43}):

\noindent
$[\bk:\bq]\le 25$ for $s=2$, $k=r=3$;

\noindent
$[\bk:\bq]\le 16$ for $s=2$, $k=3$, $r=4$;

\noindent
$\bff_5=\bq(\sqrt{5})\subset \bk$ and 
$[\bk:\bff_5]\le 8$ for $s=2$, $k=3$, $r=5$;

\noindent
$[\bk:\bq]\le 16$ for $s=2$, $k=4$, $r=3$;

\noindent
$\bff_5=\bq(\sqrt{5})\subset \bk$ and 
$[\bk:\bff_5]\le 10$ for $s=2$, $k=5$, $r=3$; 

\noindent
$[\bk:\bq]\le 17$ for $s=k=r= 3$; 

\noindent
$[\bk:\bq]\le 8$ for $s=k=3$, $r= 4$;

\noindent
$\bk=\bff_5=\bq(\sqrt{5})$ if $s=k=3$, $r=5$.

\noindent
$[\bk:\bq]\le 16$ for $s=3$, $k=4$, $r=3$;

\noindent
$\bff_5=\bq(\sqrt{5})\subset \bk$ and 
$[\bk:\bff_5]\le 10$ for $s=3$, $k=5$, $r=3$;

\noindent
$[\bk:\bq]\le 16$ for $s=k=4$, $r=3$;

\noindent
$\bff_5=\bq(\sqrt{5})\subset \bk$ and 
$[\bk:\bff_5]\le 8$ for $s=4$, $k=5$, $r=3$;

\noindent
$\bff_5=\bq(\sqrt{5})\subset \bk$ and 
$[\bk:\bff_5]\le 10$ for $s=k=5$, $r=3$;

In particular, $[\bk:\bq]\le 25$ for all fields
$\bk\in \F \Gamma_3^{(4)}(14)$.

\label{thforF43and14}
\end{theorem}

By considering the graphs $\Gamma_2^{(4)}$ directly in other 
cases, one can significantly improve Theorem \ref{propforF43and14}. 
We hope to consider this in further variants of the paper and 
further publications. 


\subsection{Fields from $\F \Gamma_5^{(4)}(14)$}\label{subsec:gam(4)5}
For $\Gamma_5^{(4)}(14)$ (see Figure \ref{graphg45}),
$3\le s\le k$ are natural numbers and $2<u<14$ is a totally real
algebraic integer.

The ground field $\bk=\bq(u^2)$ contains cyclic products
$$
\cos^2{\frac{\pi}{k}},\ \cos^2{\frac{\pi}{s}},\ u^2\ .
$$
Thus, $\bff_{s,k}=\bq(\cos^2(\pi/s),\cos^2(\pi/k))\subset \bk=\bq(u^2)$. 
The determinant $d(u)$ of the Gram matrix is determined
by the equality
$$
-\frac{d(u)}{4}=u^2-4\sin^2 {\frac{\pi}{k}}\sin^2{\frac{\pi}{s}}.
$$
The $\Gamma^{(4)}_5$ is V-arithmetic if and only if for $\alpha=u^2$ we have 
$$
0<\sigma(\alpha )<\sigma(4\,\sin^2{\frac{\pi}{k}}\sin^2{\frac{\pi}{s}})<4
$$
for any $\sigma:\bk\to \br$ which is different from identity $\sigma^{(+)}$.
For $\sigma^{(+)}$, we have
$$
4<\sigma^{(+)}(\alpha )< 14^2.
$$

\begin{figure}
\begin{center}
\includegraphics[width=6cm]{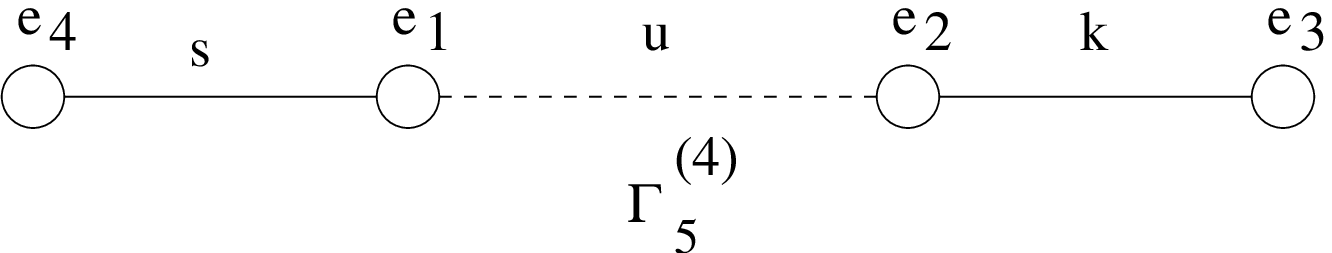}
\end{center}
\caption{The graph $\Gamma_5^{(4)}$}
\label{graphg45}
\end{figure}

The graph $\Gamma_5^{(4)}(14)$ with parameters $s$, $k$, 
contains two subgraphs $\Gamma_2^{(3)}(14)$, with parameters $d=s$ and $d=k$. 
Thus the field $\bk \in \F \Gamma_5^{(4)}(14)$ with parameters $3\le s \le k$ 
belongs to $\F \Gamma_2^{(3)}(14)$ with the parameter $s$ and to 
$\F \Gamma_2^{(3)}(14)$ with the parameter $k$. By Theorem 
\ref{thforF32and14}, the degree $[\bk:\bq]\le 20$ if one of $s$ and $k$ 
does not belong to $\{3,\,7\}$.  
Let us assume that both $s$ and $k$ belong to $\{3,\,7\}$. 

\medskip

Let $s=k=3$. Then $0<\sigma(\alpha)<9/4$.
We take  the polynomial
\begin{equation}
P(x)=x(x-1)^2(x-2)^2(x^2-3x+1).
\label{pol54for3,3}
\end{equation}
The maximum of $|P(x)|$ on the interval $[0,9/4]$ is equal to
$\delta=0.21570956...$. By Lemma \ref{lempolmain}
$$
[\bk:\bq]\le 1+\frac{\ln{P(14^2)}}{-\ln{\delta}}\le 25.06 
$$
and $[\bk:\bq]\le 25$. 

\medskip 

Let $s=3$, $k=7$. In this case,
$$
0<\sigma(\alpha)<3\sigma(\sin^2{(\pi/7)})
$$
for $\sigma:\bff_7=\bq(\cos^2(\pi/7))\to \br$.
We take the polynomial
\begin{equation}
P(x)=x^2(x-1)^2(x-2)(x^2-3x+1)(x+1-4\sin^2{(\pi/7)})(x-4\sin^2{(\pi/7)})
\label{pol54for3,7}
\end{equation}
which can be considered as an appropriate combination of the polynomials 
\eqref{pol3} for $d=3$ and \eqref{polynPd} for $d=7$. 
For $\sigma=\sigma_k:\bff_7\to \br$
where $\sigma_k(\sin^2{(\pi/7)})=\sin^2{(k\pi/7)}$, $k=1,2,3$,
we obtain that the maximum of $|P^\sigma(x)|$ on the interval
$[0,3\sigma_k(\sin^2{(\pi/7)})]$ is
equal to $\delta_1=0.0050709048...$, $\delta_2=0.110711589...$,
$\delta_3=1.23817314...$ respectively. By Lemma \ref{lempolmain}, 
\begin{equation}
[\bk:\bff_7]\le 
\frac{\ln{(P(14^2))}-\ln{(\delta_1)}}{-\ln{(\delta_1\delta_2\delta_3)}}\le 7.3.
\label{bound14for3,7}
\end{equation}
Thus $[\bk:\bff_7]\le 7$ and $[\bk:\bq]\le 21$.

Case $s=k=7$. In this case, $0<\sigma(\alpha)<\sigma(4\sin^4{(\pi/7)})$. We
take the polynomial
\begin{equation}
P(x)=x(x-1)(x+1-4\sin^2{(\pi/7)})(x-4\sin^2{(\pi/7)})
\label{pol54for7,7}
\end{equation}
which is the same as we used in \eqref{polynPd} for  $d=7$. 
For $\sigma=\sigma_k$, $k=1,2,3$, in notations
above, $|P^\sigma(x)|$ on the interval $[0,\sigma(4\sin^4{(\pi/7)})]$
has maximum $\delta_1=0.0289099468...$, $\delta_2=0.52203650...$, 
$\delta_3=2.93338845...$.
Like above, we get $[\bk:\bff_7]\le 7$ and $[\bk:\bq]\le 21$.

Thus, we obtain the following result. 

\begin{theorem}
 We have the following upper bounds for the degrees
of ground fields $\bk\in \F \Gamma_5^{(4)}(14)$ that is for
V-arithmetic edge graphs
  $\Gamma_5^{(4)}$ of the minimality $14$ (i.e., $2<u<14$) depending
on their parameters $s,k$ (see Figure \ref{graphg45}):

\noindent
$[\bk:\bq]\le 25$ for $s=k=3$;

\noindent
$[\bk:\bq]\le 21$ for $s=3$, $k=7$;

\noindent
$[\bk:\bq]\le 21$ for $s=k=7$

\noindent
$[\bk:\bq]\le 20$ 
for all other $3\le s\le k$. 

In particular, $[\bk:\bq]\le 25$ for all fields
$\bk\in \F \Gamma_5^{(4)}(14)$.

\label{thforF45and14}
\end{theorem}

This Theorem can be significantly improved using our results and methods 
from \cite{Nik6}---\cite{Nik8}. We hope to do that in further 
variants of this paper and further publications.

\subsection{Fields from  $\F \Gamma_4^{(4)}(14)$}\label{subsec:gam(4)4}
For $\Gamma_4^{(4)}(14)$ (see Figure \ref{graphg44}),
$k\ge 2$, $s,\,r \ge 3$ are natural numbers and $2<u<14$ is
a totally real algebraic integer. Moreover, we have only the following
possibilities: $s=3$, $r=3,\,4,\,5$; $s=4,5$, $r=3$.

\begin{figure}
\begin{center}
\includegraphics[width=6cm]{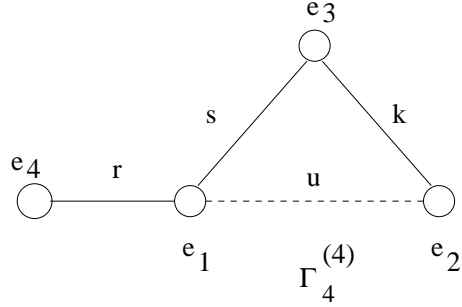}
\end{center}
\caption{The graph $\Gamma_4^{(4)}$}
\label{graphg44}
\end{figure}

The ground field $\bk=\bq(u^2)$ contains cyclic products
$$
\cos^2{\frac{\pi}{s}},\ \cos^2{\frac{\pi}{r}},\ \cos^2{\frac{\pi}{k}},\
u^2,\ u\,\cos{\frac{\pi}{s}}\cos{\frac{\pi}{k}}\,.
$$

This case can be considered as a specialization of the graph
$\Gamma^{(4)}_1$ when we take $p=2$. The determinant $d(u)$ of the
Gram matrix is determined by the equality
$$
\frac{-d(u)}{4}=u^2+4\cos{\frac{\pi}{s}}\cos{\frac{\pi}{k}}\,u+4\,
\cos^2{\frac{\pi}{s}}-4\sin^2 {\frac{\pi}{k}}\sin^2{\frac{\pi}{r}}.
$$
We obtain
\begin{equation}
\sigma(u^2+4\cos{\frac{\pi}{s}}\cos{\frac{\pi}{k}}\,u+4\,
\cos^2{\frac{\pi}{s}}-4\sin^2 {\frac{\pi}{k}}\sin^2{\frac{\pi}{r}})<0
\label{condforGamma4(4)(14)1}
\end{equation}
for $\sigma\not=\sigma^{(+)}$, and
\begin{equation}
4<\sigma^{(+)}(u^2)<14^2
\label{condforGamma4(4)(14)2}.
\end{equation}
We then have $\bk=\bq(u^2)$.

Assume that $k=2$. Then we obtain
$$
0<\sigma(u^2)<4-\sigma(4\cos^2{\frac{\pi}{s}})-\sigma(4\cos^2{\frac{\pi}{r}}).
$$
If $s=r=3$, we obtain $0<\sigma(u^2)<2$.
Applying (like above) the polynomial \eqref{pol4},
$P(x)=x(x-1)^2(x-2)$, we obtain $[\bk:\bq]\le 16$.

If $(s,r)=(3,4)$, $(4,3)$, we obtain $0<\sigma(u^2)<1$. Applying
the polynomial $P(x)=x(x-1)$, we get $[\bk:\bq]\le 8$.

If $(s,r)=(3,5)$ or $(5,3)$, we obtain
$0<\sigma(u^2)<\sigma(4\sin^2(\pi/5))-1$. Applying the polynomial
$P(x)=x(x-1)(x+1-4\sin^2{(\pi/5)})$, we obtain
$[\bk:\bff_5]\le 5$ and $[\bk:\bq]\le 10$.

Assume that $k\ge 3$. We have that
\begin{equation}
\widetilde{u}=u^2+4\cos{\frac{\pi}{s}}\cos{\frac{\pi}{k}}\,
u+4\cos^2{\frac{\pi}{s}}
\label{condforGamma4(4)(14)3}
\end{equation}
is a totally positive algebraic integer
(since minimum of this quadratic polynomial of $u$ is equal to
$4\cos^2{\frac{\pi}{s}}\sin^2{\frac{\pi}{k}}$) which belongs to $\bk$,
and $\Gamma^{(4)}_4$
is V-arithmetic if and only if
\begin{equation}
0<\sigma(4\cos^2{\frac{\pi}{s}}\sin^2{\frac{\pi}{k}})\le
\sigma(\widetilde{u})<\sigma(4\,\sin^2{\frac{\pi}{r}}\sin^2{\frac{\pi}{k}}) <4
\label{condforGamma4(4)(14)4}
\end{equation}
for any $\sigma:\bk\to \br$ which is different from identity $\sigma^{(+)}$.
For $\sigma^{(+)}$, we have
\begin{equation}
4<4+8\cos{\frac{\pi}{s}}\cos{\frac{\pi}{k}}+4\cos^2{\frac{\pi}{s}}
<\sigma^{(+)}(\widetilde{u})<
14^2+56\cos{\frac{\pi}{s}}\cos{\frac{\pi}{k}}+4\cos^2{\frac{\pi}{s}}
<16^2.
\label{condforGamma4(4)(14)5}
\end{equation}
It follows that $\bk=\bq(\widetilde{u})$.

Further we use only the inequalities
$$
0<\sigma(\widetilde{u})<\sigma(4\,\sin^2{\frac{\pi}{r}}\sin^2{\frac{\pi}{k}})
<4
$$
and
$$
4<\sigma^{(+)}(\widetilde{u})<16^2.
$$
Then this case is similar to the previous case of $\Gamma^{(4)}_5(14)$
if we replace $14^2$ by $16^2$
(or use \eqref{condforGamma4(4)(14)5} if necessary),
and take $s=\min{(r,k)}$.
Instead of Theorem \ref{thforF32and14} we should use 
Theorem \ref{thforF32and16}. 

Here similar bad cases are as follows:

Case $r=k=3$. Then $[\bk:\bq]\le 25$ (here we use
$14^2+56\cos{\frac{\pi}{s}}\cos{\frac{\pi}{k}}+4\cos^2{\frac{\pi}{s}}$
instead of $14^2$).

Case $r=3$, $k=7$. Then $[\bk:\bq]\le 21$

Thus, we obtain the following.

\begin{theorem}
 We have the following upper bounds for the degrees
of ground fields $\bk\in \F \Gamma_4^{(4)}(14)$ that is for
V-arithmetic edge graphs
  $\Gamma_4^{(4)}$ of the minimality $14$ (i.e., $2<u<14$) depending
on their parameters $k,s,r$ (see Figure \ref{graphg45}):

\noindent
$[\bk:\bq]\le 16$ for $k=2$, $s=r=3$;

\noindent
$[\bk:\bq]\le 8$ for $k=2$, $(s,r)=(3,4),\,(4,3)$;

\noindent
$[\bk:\bq]\le 10$ for $k=2$, $(s,r)=(3,5),\,(5,3)$;

\noindent
$[\bk:\bq]\le 25$ for $r=k=3$;

\noindent
$[\bk:\bq]\le 21$ for $r=3$, $k=7$;

\noindent
$[\bk:\bq]\le 20$ 
for all remaining $k,s,r$.  

In particular, $[\bk:\bq]\le 25$ for all fields
$\bk\in \F \Gamma_4^{(4)}(14)$.

\label{thforF44and14}
\end{theorem}

This Theorem can be significantly improved using our results and methods 
from \cite{Nik6}---\cite{Nik8}. We hope to demonstrate that in further 
variants of this paper and further publications. 

\medskip 

This finishes the proof of Theorem \ref{thdegrees}.

\section{Applications}
\label{Applications}

Using results of \cite{Nik1}---\cite{Nik8}, we obtain the 
following applications to arithmetic hyperbolic reflection groups.

\begin{theorem} In dimensions $n\ge 10$, the ground field of
any arithmetic hyperbolic reflection group belongs to one of finite
sets of fields $\F L^4$, $\F T$, $\F\Gamma^{(4)}_i(14)$, $1\le i\le 5$. 
In particular, its degree has the upper
bound $25$.
\label{thfor10}
\end{theorem}

\begin{proof} 
In \cite{Nik1}, \cite{Nik2} (see more details in \cite{Nik6}) 
we showed that the field belongs to $\F L^4$, $\F T$, 
$\F\Gamma^{(4)}_i(14)$, $1\le i\le 5$. 
By Theorem \ref{thdegrees} of this paper, than the degree 
has the upper bound $25$. (See more details in the proof of the 
next Theorem.)
\end{proof}

\begin{theorem} In dimensions $n\ge 6$, the ground field of
any arithmetic hyperbolic reflection group belongs to one of finite
sets of fields $\F L^4$, $\F T$, $\F\Gamma^{(4)}_i(14)$, $1\le i\le 5$,
and $\F\Gamma_{2,4}(14)$. In particular, its degree has the upper
bound $25$.
\label{thfor6}
\end{theorem}

\begin{proof} This is similar to \cite[Theorem 9]{Nik6} where we
claimed this theorem for $56$.

In \cite{Nik6}, we showed that the
field belongs to one of sets $\F L^4$, $\F T$,
$\F\Gamma^{(4)}_i(14)$, $1\le i\le 5$, and $\F\Gamma_{2,4}(14)$.
It is well-known that the degree of fields from $\F L^4$ has
the upper bound $2$. By Takeuchi \cite{Tak2} (see also
\cite{Tak1}---\cite{Tak4}), the degree of fields from $\F T$
has the upper bound $5$. In \cite{Nik6}, we showed
that the degree of fields from $\F \Gamma_{2,4}(14)$ has the upper bound
$11$. By Theorem \ref{thdegrees} of this paper, the degree of fields from
$\F\Gamma^{(4)}_i(14)$, $1\le i\le 5$, has the upper bound $25$.

It follows the statement.
\end{proof}

\begin{theorem} In dimensions $n=4,\, 5$, the ground field of
any arithmetic hyperbolic reflection group belongs to one of finite
sets of fields
$\F \Gamma_1^{(6)}(14)$, $\F \Gamma_2^{(6)}(14)$, $\F \Gamma_3^{(6)}(14)$,
$\F \Gamma_1^{(7)}(14)$, $\F \Gamma_2^{(7)}(14)$, and
$\F \Gamma_{2,5}(14)$ if it is different from a field from the sets
$\F L^4$, $\F T$, $\F\Gamma^{(4)}_i(14)$, $1\le i\le 5$, $\F\Gamma_{2,4}(14)$
in Theorem \ref{thfor6} above.

In particular, its degree has the upper
bound $44$.
\label{thfor4,5}
\end{theorem}

\begin{proof} This is similar to \cite[Theorem 2.6]{Nik7} where we
claimed the upper bound $138$.

In \cite{Nik7}, we showed that if a
field is different from a field from  $\F L^4$, $\F T$,
$\F\Gamma^{(4)}_i(14)$, $1\le i\le 5$, $\F\Gamma_{2,4}(14)$,
then it belongs to
$\F \Gamma_1^{(6)}(14)$, $\F \Gamma_2^{(6)}(14)$, $\F \Gamma_3^{(6)}(14)$,
$\F \Gamma_1^{(7)}(14)$, $\F \Gamma_2^{(7)}(14)$, and
$\F \Gamma_{2,5}(14)$. In \cite{Nik7} we showed that the degree of
fields from $\F \Gamma_{2,5}(14)$ has the upper bound $12$.

Assume that a field belongs to
$\F \Gamma_1^{(6)}(14)$, $\F \Gamma_2^{(6)}(14)$, $\F \Gamma_3^{(6)}(14)$,
$\F \Gamma_1^{(7)}(14)$, $\F \Gamma_2^{(7)}(14)$. 
All diagrams $\Gamma_1^{(6)}(14)$, $\Gamma_2^{(6)}(14)$, $\Gamma_3^{(6)}(14)$,
$\Gamma_1^{(7)}(14)$, $\Gamma_2^{(7)}(14)$ have a subdiagram
$\Gamma_2^{(3)}(14)$. By Theorem \ref{thfieldsGamma2(3)}, then the
degree of the field has the upper bound $44$.

It follows the statement.
\end{proof}

\begin{theorem} In dimension $n=3$, the ground field of
any arithmetic hyperbolic reflection group belongs to one of finite
sets of fields
$\F \Gamma_6^{(4)}(14)$, \linebreak 
$\F \Gamma_1^{(5)}(14)$, $\F \Gamma_4^{(6)}(14)$,
and $\F \Gamma_2(14)$ if it is different from a field from
the sets
$\F \Gamma_1^{(6)}(14)$, $\F \Gamma_2^{(6)}(14)$, $\F \Gamma_3^{(6)}(14)$,
$\F \Gamma_1^{(7)}(14)$, $\F \Gamma_2^{(7)}(14)$,
$\F \Gamma_{2,5}(14)$ and
$\F L^4$, $\F T$, $\F\Gamma^{(4)}_i(14)$, $1\le i\le 5$, $\F\Gamma_{2,4}(14)$
in Theorems \ref{thfor6} and \ref{thfor4,5} above.

In particular, its degree has the upper
bound $44$.
\label{thfor3}
\end{theorem}

\begin{proof} This is similar to \cite[Theorem 3.8]{Nik8} where we
claimed the upper bound $909$.

In \cite{Nik8}, we showed that if a
field is different from a field from  $\F L^4$, $\F T$,
$\F\Gamma^{(4)}_i(14)$, $1\le i\le 5$, $\F\Gamma_{2,4}(14)$, and
$\F \Gamma_1^{(6)}(14)$, $\F \Gamma_2^{(6)}(14)$, $\F \Gamma_3^{(6)}(14)$, 
\linebreak 
$\F \Gamma_1^{(7)}(14)$, $\F \Gamma_2^{(7)}(14)$, $\F \Gamma_{2,5}(14)$,
then it belongs to one of sets
$\F \Gamma_6^{(4)}(14)$, $\F \Gamma_1^{(5)}(14)$, $\F \Gamma_4^{(6)}(14)$,
and $\F \Gamma_2(14)$.

Assume that a field belongs to one of sets
$\F \Gamma_6^{(4)}(14)$, $\F \Gamma_1^{(5)}(14)$, $\F \Gamma_4^{(6)}(14)$,
and $\F \Gamma_2(14)$.

The set $\F \Gamma_2(14)$ consists of ground fields of arithmetic hyperbolic
reflection groups in dimension 2 with the fundamental polygon of 
the minimality $14$. We showed in 
\cite{Nik6} using results of \cite{LMR} that 
the degree of ground fields of arithmetic hyperbolic 
reflection groups of dimension $2$ has the upper bound $44$. 

All diagrams
$\Gamma_6^{(4)}(14)$, $\Gamma_1^{(5)}(14)$, $\Gamma_4^{(6)}(14)$
have a subdiagram
$\Gamma_2^{(3)}(14)$. By Theorem \ref{thfieldsGamma2(3)}, then the
degree of the field has the upper bound $44$.

It follows the statement.
\end{proof}

We remark that Theorem \ref{thfor3} also gives another proof of
finiteness in dimension $3$ which is different from the proof by I. Agol
\cite{Agol}.

In \cite{B}, Belolipetsky showed (using results of \cite{Agol}) that
the degree of ground fields of arithmetic hyperbolic reflection groups
in dimension $3$ has the upper bound $35$. Using this result and
results of \cite{Nik6},
we can improve the upper bound above in dimensions $4$ and $5$.

\begin{theorem} In dimensions $n=4$, $5$, the ground field of
an arithmetic hyperbolic reflection group belongs to one of
sets $\F L^4$, $\F T$, $\F\Gamma^{(4)}_i(14)$, $1\le i\le 5$,
or it is equal to the ground field of an
arithmetic hyperbolic reflection group in dimension $2$ or $3$.

In particular, its degree has the upper bound $35$.
\label{thforalldim}
\end{theorem}

\begin{proof} In \cite{Nik5} we proved the first statement of
the theorem.

The degree of fields from $\F L^4$, $\F T$ has the upper bound $5$.
By Theorem \ref{thdegrees} of this paper, the degree of fields from
$\F\Gamma^{(4)}_i(14)$, $1\le i\le 5$, has the upper bound $25$.

By \cite{B}, the degree of fields in dimension $3$ has the
upper bound $35$.

By \cite{M}, the degree of fields in dimension $2$ has the upper
bound $11$.

It follows the statement.
\end{proof}

Theorems \ref{thfieldsGamma2(3)} and \ref{thforF32and14}
can be also considered as a small step to classification.

\begin{theorem} Assume that a narrow face of the minimality $14$
of the  fundamental chamber (it does exist by \cite{Nik1}) 
of an arithmetic hyperbolic reflection group contains an edge with 
the hyperbolic connected component of its graph containing 
a subgraph $\Gamma_2^{(3)}$ with the parameter
$d$ (see Figure \ref{graphg3}).

Then $d\le 120$, and the ground field is $\bff_d=\bq(\cos^2{(\pi/d)})$
if $d\ge 44$ (see more exact results in Theorem \ref{thforF32and14}). 
\label{class1}
\end{theorem}

\section{Appendix: Hyperbolic numbers (the review of \cite[Sec. 1]{Nik2})}
\label{sec:Appendix}

Here we review our results in \cite[Sec.1]{Nik2} and correct some
arithmetic mistakes (Theorems 1.1.1 and 1.2.2 in \cite{Nik2}
which are similar
to Theorems \ref{thFekete2} and \ref{thhypnumberscond} here).
This mistakes are unessential for results of \cite{Nik2}.

\subsection{Fekete's theorem}\label{subsec:fekete}
The following
important theorem, to which this section is devoted, was obtained
by Fekete \cite{Fek}. Although Fekete considered (as we know) the
case of $\bq$ his method of proof can be immediately carried over
to totally real algebraic number fields.

\begin{theorem} (M. Fekete). Suppose that $\bff$ is a totally real
algebraic number field, and to every embedding $\sigma: \bff\to
\br$ there corresponds an interval $[a_\sigma,b_\sigma]$ in $\br$
and the real number $\lambda_\sigma>0$. Suppose that
$\prod_\sigma{\lambda_\sigma}=1$. Then for every nonnegative
integer $n$ there exists a nonzero polynomial $P_n(T)\in \bo[T]$
of degree no greater than $n$ over the ring of integers $\bo$ of
$\bff$ such that the following inequality holds for each $\sigma$:
\begin{equation}
\max_{x\in [a_\sigma,b_\sigma]}{|P_n^\sigma(x)|}\le \lambda_\sigma
 {|\discr \bff|^{1/(2[\bff:\bq])}
2^{n/(n+1)}(n+1)\left(\prod_\sigma {\frac{b_\sigma -
a_\sigma}{4}}\right)^{n/(2[\bff:\bq])}}
\label{fek1a}\ .
\end{equation}
\label{thFekete1}
\end{theorem}

\begin{proof} Suppose that $N=[\bff:\bq]$ and that
$\gamma_1,\dots\gamma_N$ is the basis for $\bo$ over $\bz$.
Suppose we are given a nonzero polynomial
$$
P_n(T)=\sum_{i=0}^{n}{\sum_{j=1}^{N}{\alpha_{ij}\gamma_jT^i}}\in
\bo[T]
$$
of degree no greater than $n$, where the $\alpha_{ij}\in \bz$
are not all zero. For every $\sigma:\bff\to \br$ we consider the
real functions $P_n^\sigma(x)$ on the interval
$[a_\sigma,b_\sigma]$.

We make the change of variables
$$
x=x(z)=\frac{b_\sigma+a_\sigma}{2}+\frac{b_\sigma-a_\sigma}{2}\cos{z}.
$$
If $z$ runs through $[0,\pi]$, then $x$ runs through
$[a_\sigma,b_\sigma]$. We also set
$Q^\sigma_n(z)=P_n^\sigma(x(z))$.

Since $Q_n^\sigma(z)$ is an even trigonometric polynomial, it
follows that
\begin{equation}
Q_n^\sigma(z)=\sum_{k=0}^{n}A_{k\sigma}\cos{kz},
\label{Qnsigma}
\end{equation}
where
$$
A_{k\sigma}=\frac{1}{\pi}\int_{-\pi}^{\pi}
{P_n^\sigma\left(\frac{b_\sigma+a_\sigma}{2}+\frac{b_\sigma-a_\sigma}{2}\cos{z}
\right)\cos{kz}\,dz}
$$
$$
=\sum_{i=0}^{n}{\sum_{j=1}^{N}\left(\frac{1}{\pi}\int_{-\pi}^{\pi}{\gamma_j^\sigma
\left(\frac{b_\sigma+a_\sigma}{2}+\frac{b_\sigma-a_\sigma}{2}\cos{z}
\right)^i\cos{kz}\,dz}\right)}\alpha_{ij}\,,
$$
if $k\ge 1$, and
$$
A_{0\sigma}=\frac{1}{2\pi}\int_{-\pi}^{\pi}
{P_n^\sigma\left(\frac{b_\sigma+a_\sigma}{2}+\frac{b_\sigma-a_\sigma}{2}\cos{z}
\right)\,dz}
$$
$$
=\sum_{i=0}^{n}{\sum_{j=1}^{N}\left(\frac{1}{2\pi}\int_{-\pi}^{\pi}
{\gamma_j^\sigma
\left(\frac{b_\sigma+a_\sigma}{2}+\frac{b_\sigma-a_\sigma}{2}\cos{z}
\right)^i \,dz}\right)}\alpha_{ij}\,.
$$
Thus,
\begin{equation}
A_{k\sigma}=\sum_{i=0}^{n}{\sum_{j=1}^{N}{c_{k\sigma
ij}\alpha_{ij}} \label{Aksigma}}
\end{equation}
are linear functions of the $\alpha_{ij}$, where
$$
c_{k\sigma ij}=\gamma_j^\sigma\cdot \frac{1}{\pi}\int_{-\pi}^{\pi}
{\left(\frac{b_\sigma+a_\sigma}{2}+\frac{b_\sigma-a_\sigma}{2}\cos{z}
\right)^i\cos{kz}\,dz}\, ,
$$
if $k\ge 1$, and
$$
c_{0\sigma ij}=\gamma_j^\sigma\cdot
\frac{1}{2\pi}\int_{-\pi}^{\pi}
{\left(\frac{b_\sigma+a_\sigma}{2}+\frac{b_\sigma-a_\sigma}{2}\cos{z}
\right)^i \,dz}\, .
$$
We note that, because of these formulas, $c_{k\sigma ij}=0$ for
$i<k$, and
$$
c_{k\sigma kj}=\gamma_j^\sigma\cdot
2\left(\frac{b_\sigma-a_\sigma}{4}\right)^k \,, \text{\ \  if
}k\ge 1
$$
($c_{0\sigma0j}=\gamma_j^\sigma$). Hence, if we order the indices
$k\sigma$ and $ij$ lexicographically, we find that the matrix of
the linear forms \eqref{Aksigma} is an upper block-triangular
matrix with the shown above $N\times N$ matrices $(c_{0\sigma 0 j})$ and
$(c_{k\sigma k j})$, $1\le k \le n$, on
the diagonal. It follows that its determinant is equal to
$$
\Delta=\det(\gamma_j^{\sigma})^{n+1}\cdot
2^{Nn}\left(\prod_\sigma{\frac{b_\sigma-a_\sigma}{4}}\right)^{n(n+1)/2}\,.
$$
Since $\prod_\sigma \lambda_\sigma^{n+1}=
\left(\prod_\sigma \lambda_\sigma  \right)^{n+1}=1$,
according to Minkowski's theorem on linear forms (see, for
example, \cite{BS}, \cite{Cas}), there exist $\alpha_{ij}\in \bz$, not
all zero, such that
$|A_{k\sigma}|\le
\lambda_\sigma |\Delta|^{1/N(n+1)}$,
and hence, by \eqref{Qnsigma},
$$
\max_{z}{\vert Q_n^\sigma(z)\vert }\le
\lambda_\sigma \cdot (n+1)|\Delta|^{1/N(n+1)}.
$$
Taking into account that $\det(\gamma_j^\sigma)^2=\discr \bff$, we
obtain the proof of the theorem.
\end{proof}

Taking $\lambda_\sigma=1$, we get a particular statement which we later use.

\begin{theorem} (M. Fekete). Suppose that $\bff$ is a totally real
algebraic number field, and to every embedding $\sigma: \bff\to
\br$ there corresponds an interval $[a_\sigma,b_\sigma]$ in $\br$.
Then for every nonnegative
integer $n$ there exists a nonzero polynomial $P_n(T)\in \bo[T]$
of degree no greater than $n$ over the ring of integers $\bo$ of
$\bff$ such that the following inequality holds for each $\sigma$:
\begin{equation}
\max_{x\in [a_\sigma,b_\sigma]}{|P_n^\sigma(x)|}\le
 {|\discr \bff|^{1/(2[\bff:\bq])}
2^{n/(n+1)}(n+1)\left(\prod_\sigma {\frac{b_\sigma -
a_\sigma}{4}}\right)^{n/(2[\bff:\bq])}}
\label{fek1}\ .
\end{equation}
\label{thFekete2}
\end{theorem}

\subsection{Hyperbolic numbers}
\label{hypnumbers}
The totally real algebraic integers $\{\alpha\}$ which we consider
here are very similar to Pisot-Vijayaraghavan numbers \cite{Cas}, although
the later are not totally real.

\begin{theorem} Let $\bff$ be a totally real algebraic number field,
and let each embedding $\sigma:\bff \to \br$ corresponds
to an interval $[a_\sigma,b_\sigma]$ in $\br$, where
$$
\prod_\sigma{\frac{b_\sigma-a_\sigma}{4}}<1.
$$
In addition, let the natural number $m$ and the intervals
$[s_1,t_1],\dots , [s_m,t_m]$ in $\br$
be fixed.

Then there exists a constant $N(s_i,t_i)$ such that, if $\alpha$
is a totally real algebraic integer
and if the following inequalities hold for the embeddings
$\tau:\bff(\alpha)\to \br$:
$$
s_i\le \tau(\alpha)\le t_i,\ \text{for\ }\tau=\tau_1,\dots, \tau_m,
$$
$$
a_{\tau\vert \bff}\le \tau(\alpha)\le b_{\tau\vert \bff}\
\text{for\ }\tau\not=\tau_1,\dots\tau_m,
$$
then
$$
[\bff(\alpha):\bff]\le N(s_i,t_i).
$$
\label{thhypnumbers}
\end{theorem}

\begin{theorem}
Under the conditions of Theorem \ref{thhypnumbers}, $N(s_i,t_i)$
can be taken to be $N(s_i,t_i)=N(S)$, where $N(S)$ is the least natural
number solution of the inequality
\begin{equation}
nln(1/R)-M\ln(2n+2)-\ln{B}\ge \ln{S}.
\label{ineqforn}
\end{equation}
Here
\begin{equation}
M=[\bff:\bq],\ \ \
R=\sqrt{\prod_\sigma{\frac{b_\sigma-a_\sigma}{4}}},
\label{ineqforMR}
\end{equation}
\begin{equation}
B=\sqrt{\vert \discr \bff\vert},\ \ \
S=\prod_{i=1}^{m}{\left( 2e r_i(b_{\sigma_i}-a_{\sigma_i})^{-1}\right)},
\label{ineqforBS}
\end{equation}
where
$\sigma_i=\tau_i\vert\bff$ and $r_i=\max{\{\vert t_i-a_{\sigma_i}\vert,\
\vert b_{\sigma_i}-s_i\vert\}}$.

Asymptotically,
$$
N(s_i, t_i)\sim \frac{\ln S }{\ln {(1/R)}}.
$$
\label{thhypnumberscond}
\end{theorem}

\begin{proof} We use the following statement.

\begin{lemma} Suppose that $Q_n(T)\in \br[T]$ is a non-zero polynomial
over $\br$ of degree no greater than $n>0$,
$a<b$ and $M_0=\max_{[a,b]}{\{\vert Q_n(x)\vert\}}$. Then for $x\ge b$
$$
|Q_n(x)|\le \frac{M_0(x-a)^n n^n }{\left((b-a)/2\right)^n n!}<
$$
$$
\frac{M_0(x-a)^n e^n}{((b-a)/2)^n\sqrt{2\pi n}} <
\frac{M_0(x-a)^n e^n}{((b-a)/2)^n} .
$$
\label{lemLagr}
\end{lemma}

\begin{proof} Let $\alpha_0<\alpha_1<\cdots <\alpha_n$.
Then we have the Lagrange interpolation formula
$$
Q_n(x)=\sum_{i=0}^{n}{Q_n(\alpha_i)F_i(x)}
$$
where
$$
F_i(x)=\frac{(x-\alpha_0)(x-\alpha_1)\cdots
(x-\alpha_{i-1})(x-\alpha_{i+1})\cdots (x-\alpha_n)}
{(\alpha_i-\alpha_0)(\alpha_i-\alpha_1)\cdots
(\alpha_i-\alpha_{i-1})(\alpha_i-\alpha_{i+1})\cdots(\alpha_i-\alpha_n)}.
$$
Taking $\alpha_i=a+i(b-a)/n$, $0\le i\le n$, we obtain for $x\ge b$ that
$$
|Q_n(x)|\le \frac{M_0(x-a)^n}{((b-a)/n)^n}\sum_{i=0}^{n}\frac{1}{i!(n-i)!}=
\frac{M_0(x-a)^n2^n}{((b-a)/n)^n n!}.
$$
By Stirling formula, $n!=\sqrt{2\pi n}(n/e)^ne^{\lambda_n}$ where
$0<\lambda_n<1/(12n)$. Thus, $n^n/n!< e^n/\sqrt{2\pi n}<e^n$. It
follows the statement.
\end{proof}

We continue the proof of theorems.

For given $n$ we consider the polynomial $P_n(T)\in \bo[T]$ whose
existence is ensured by Fekete's theorem \ref{thFekete2}. Setting
$N=[\bff(\alpha):\bff]$ and $M=[\bff:\bq]$, we use Fekete's
theorem and the lemma to conclude that
$$
|\prod_\tau{\tau(P_n(\alpha))|}=\prod_{\tau}{|P_n^\tau(\tau(\alpha))|}=
\prod_{\tau\not=\tau_i}{P_n^\tau(\tau(\alpha))}
\prod_{i=1}^{m}{|P_n^{\tau_i}(\tau_i(\alpha))|}
$$
$$
\le \prod_{\tau\not=\tau_i}{\max_{[a_{\tau|\bff},b_{\tau|\bff}]}{|P_n^\tau(x)|}}
\prod_{i=1}^m{\max_{[s_i,t_i]}{|P_n^{\sigma_i}(x)|}}
$$
$$
\le\left(|\discr \bff|^{1/(2M)}\cdot 2\cdot (n+1)R^{n/M}\right)^{NM}
\prod_{i=1}^{m}{\frac{r_i^ne^n}{((b_{\sigma_i}-a_{\sigma_i})/2)^n}}
$$
$$
=R^{nN}B^N\cdot S^n\cdot (2n+2)^{MN}.
$$
Since $R<1$, there exists $n_0$ large enough so that
\begin{equation}
R^{n_0}\cdot B\cdot(2n_0+2)^{M}\le \frac{1}{S}.
\label{eqn0}
\end{equation}
Then if $N>n_0$, we find that
$$
R^{n_0N}\cdot B^N\cdot S^{n_0}(2n_0+2)^{MN}\le S^{n_0-N}<1,
$$
since $S>1$. From this and the above chain of inequalities we have
$$
|\prod_\tau{\tau(P_{n_0}(\alpha))}|<1.
$$
But
$$
\prod_\tau{\tau(P_{n_0}(\alpha))}=N_{\bff(\alpha)/\bq}(P_{n_0}(\alpha))\in \bz,
$$
and hence $P_{n_0}(\alpha)=0$. Consequently, $N\le n_0$, and we have
obtained a contradiction.
We have thereby proved that $N\le n_0$, where $n_0$ is a natural
number solution of \eqref{eqn0}.
The inequality \eqref{eqn0} is obviously equivalent to
$$
n_0\ln{(1/R)}-M\ln{(2n_0+2)}-\ln{B}\ge \ln{S},
$$
and this completes the proof of the theorems.
\end{proof}


V.V. Nikulin \par Deptm. of Pure Mathem. The University of
Liverpool, Liverpool\par L69 3BX, UK; \vskip1pt Steklov
Mathematical Institute,\par ul. Gubkina 8, Moscow 117966, GSP-1,
Russia

vnikulin@liv.ac.uk \ \ vvnikulin@list.ru

\end{document}